\documentclass[a4paper, 11pt]{article}
\usepackage{hyperref}
\usepackage{indentfirst}
\usepackage{epsfig}
\usepackage{psfrag}
\usepackage{graphicx}
\usepackage{overpic}
\usepackage{multirow}
\usepackage[hang]{subfigure}
\usepackage{amsmath,amssymb, amsthm}
\usepackage{color}
\usepackage[mathscr]{euscript}

\newtheorem{definition}{Definition}

\newcommand\bbR{\mathbb{R}}
\newcommand\bbN{\mathbb{N}}
\newcommand\bbS{\mathbb{S}}

\newcommand\bg{\boldsymbol{g}}
\newcommand\bn{\boldsymbol{n}}
\newcommand\bu{\boldsymbol{u}}
\newcommand\bv{\boldsymbol{v}}
\newcommand\bw{\boldsymbol{w}}
\newcommand\bx{\boldsymbol{x}}
\newcommand\bD{{\bf{D}}}
\newcommand\bE{\boldsymbol{E}}
\newcommand\bF{\boldsymbol{F}}
\newcommand\bM{{\bf{M}}}
\newcommand\dd{\,\mathrm{d}}
\newcommand\He{\mathit{He}}

\newcommand\mC{{\mathcal{C}}}
\newcommand\mH{\mathcal{H}}
\newcommand\mE{\mathcal{E}}

\newcommand\pd[2]{\dfrac{\partial {#1}}{\partial {#2}}}

\newcommand\od[2]{\dfrac{\mathrm{d} {#1}}{\mathrm{d} {#2}}}

\newcommand\rC[2]{{\rm{C}}_{{#1},{#2}}}

\newcommand\comment[1]{}

\graphicspath{{images/}}
{\theoremstyle{remark} \newtheorem{remark}{Remark}}
\newtheorem{proposition}{Proposition}

\title{Filtered Hyperbolic Moment Method for the Vlasov Equation}

\author{Yana Di\thanks{LSEC, NCMIS, Academy of Mathematics and Systems
    Science, Chinese Academy of Sciences, Beijing, 100190, China, and
    School of Mathematical Sciences, University of Chinese Academy of
    Sciences, Beijing, 100049, China, email: {\tt
      yndi@lsec.cc.ac.cn}.}, ~~ Yuwei Fan\thanks{Department of
    Mathematics, Stanford University, Stanford, CA 94305, email: {\tt
      ywfan@stanford.edu}.},~~ Zhenzhong Kou\thanks{Academy of
    Mathematics and Systems Science, Chinese Academy of Sciences,
    Beijing, China, email: {\tt kouzhzh@lsec.cc.ac.cn}.},~~ Ruo
  Li\thanks{CAPT, LMAM \& School of Mathematical Sciences, Peking
    University, Beijing, China, email: {\tt rli@math.pku.edu.cn}.},~~
  Yanli Wang\thanks{College of Engineering, Peking University,
    Beijing, China, email: {\tt wang\_yanli@pku.edu.cn}.}  }

\begin{document}
\maketitle
% vim: tw=70:spell
\begin{abstract}
  In this paper, we investigate the effect of the filter for the
  hyperbolic moment equations (HME) \cite{Fan_new, Wang} of the
  Vlasov-Poisson equations and propose a novel quasi time-consistent
  filter to suppress the numerical recurrence effect. By taking
  properties of HME into consideration, the filter preserves a lot of
  physical properties of HME, including Galilean invariance and 
  conservation of mass, momentum and energy. We present two
  viewpoints --- collisional viewpoint and dissipative viewpoint --- to
  dissect the filter, and show that the filtered hyperbolic moment
  method can be treated as a solver of Vlasov equation. Numerical
  simulations of the linear Landau damping and two stream instability
  demonstrate the effectiveness of the filter in restraining
  recurrence arising from particle streaming. Both the analysis and
  the numerical results indicate that the filtered method can capture the
  evolution of the Vlasov equation, even when phase mixing and
  filamentation dominate.
  
  \vspace*{4mm}
 
\end{abstract}

\section{Introduction}
The Vlasov equation is the fundamental kinetic equation modeling of
the collisionless plasma. It describes the time evolution of the
distribution $f(t,\bx,\bv)$ of a population of charged particles
(electrons, ions) that responds to the self-consistent electromagnetic
fields. The distribution function $f(t,\bx,\bv)$ is the number density
of the particles at the time $t$ and position
$\bx \in \Omega \subset \bbR^D$ with the microscopic velocity
$\bv \in \bbR^D$ \cite{Vlasov} (for physical case $D=3$). In this
paper, we focus on the dimensionless Vlasov-Poisson equations (VP)
\begin{equation} \label{eq:vlasov}
    \pd{f}{t} + \bv \cdot \nabla_{\bx} f 
    + \bF(t,\bx) \cdot \nabla_{\bv} f = 0, 
    \quad (\bx, \bv, t) \in  \Omega \times \bbR^D \times \bbR^+, 
\end{equation}
where $\bF(t,\bx)$ is the electric force produced by the
self-consistent electric filed $\bE(t,\bx)$:
\begin{equation} \label{eq:force}
    \bF(t,\bx) = \bE(t,\bx), \quad
    \bE(t,\bx) = -\nabla_{\bx} ~\phi(t,\bx),\quad
    -\Delta_{\bx} \phi(t,\bx) = \rho(t,\bx) - \rho_0,
    \end{equation}
    where $\phi(t,\bx)$ is the electric potential produced by the
    particles. For the clarity of notations, we will always assume
    periodicity in $\bx$.  The density $\rho(t, \bx)$ is defined as
    $\rho(t,\bx)=\int_{\bbR^D}f(t,\bx,\bv)\dd\bv$ and $\rho_0$ is
    positive constant dependent on the problem. When the VP system is
    applied to the plasmas, the total charge neutrality condition
    $\int_{\Omega}(\int_{\bbR^D} f \dd \bv -\rho_0) \dd \bx = 0.$

Numerical methods for solving the Vlasov equation have been
extensively studied, see for instance \cite{birdsall1985plasma,
  zaki1988finiteI, heath2012discontinuous, Filbet2003Comparison} and
references therein.  The most common method is the Particle-In-Cell
(PIC) method, where the Vlasov equation is solved by following the
trajectories of a set of statistically distributed velocity points.
The method is proved to be successful due to its relative simplicity
and adaptivity \cite{eliasson2010numerical,
  camporeale2016velocity}. But as a stochastic method, the inherent
statistical noise of PIC sometimes overshadows the physical results.
Moreover, the large tail of the distribution function would lower
PIC's effectiveness.  With the exponential growth of the computing
power, the Eulerian numerical methods attract more and more
researchers' attention, and a lot of methods are developed, for
example, the continuous finite element methods \cite{zaki1988finiteI},
finite difference methods \cite{FATEMI1993209, Carrillo2003}, finite
volume methods \cite {Filbet2001Convergence}, discontinuous Galerkin method
\cite{heath2012discontinuous, Qiu2011positivity}, the semi-Lagrangian
method \cite{SL} and spectral methods \cite{UTH,
  Bourdiec2006Numerical, Wang}.  These methods discretize or
approximate the distribution function both in the spatial space and
the microscopic velocity space, and they can be used to solve the case
that the distribution function has the low-density velocity much more
accurately compared to PIC, but may be quite expensive in the high
dimension problem.

However, a challenging issue of the Eulerian numerical methods is the
filamentation \cite{Shoucri1974Numerical} for instance, the discrete
velocity method and the spectral method.  The filamentation is caused
by the oscillations with smaller and smaller wavelengths in velocity
space of the distribution function as time evolves
\cite{camporeale2016velocity}. The filamentation phenomenon is a
common issue among the Eulerian numerical methods
\cite{Filbet2003Comparison} due to its limitation on the
representation of high frequency information.  A classical and well
studied example of filamentation is the linear Landau damping
\cite{armstrong1970solution}. The failure on capturing filamentation
leads to the well-known numerical recurrence phenomenon
\cite{Bourdiec2006Numerical}.

Even though the Eulerian numerical methods can not capture the
filamentation, but it is possible to suppress or even eliminate the
recurrence phenomenon.  The key idea to suppress the recurrence is to
capture the main structures of the particles when discarding the
information of oscillations with small wave length.  Artificial
collisional term is a popular method to suppress the recurrence
\cite{armstrong1967, Joyce1971Numerical}.  The additional collisional
term can be interpreted as an advection-diffusion operator, thus it
can eliminate the high frequency of oscillations in the velocity
space.  Several kinds of collisional terms have been adapted, for
example, the weakly Fokker-Planck collision operator was suggested in
\cite{grant1967fourier}, and a numerical collision term by a nonlinear
combination of Lenard-Bernstein collision operator was tested in
\cite{camporeale2016velocity}.  Filtering is a common procedure to
reduce the effects of the Gibbs phenomenon in spectral methods
\cite{canuto2012spectral} and is also a popular method to suppress the
recurrence \cite{Cheng, Klimas1987method, Klimas1994splitting,
  Holloway1996Spectral}.  The authors in \cite{Cheng} pointed out that
a filter of high quality could weaken the high oscillation of the
particles and have little influence on the lower order moments of the
distribution function.  In \cite{parker2015fourier}, Hou-Li filter
which was proposed in \cite{HouLi2007} was adapted in velocity space
to a Fourier-Hermite spectral representation to study Landau damping.
Moreover, some other methods, for example, the absorbing boundary
conditions \cite{ellasson2001outflow}, were proposed to retain the
recurrence.

In this paper, we focus on how to apply the filter onto the hyperbolic
moment method for the Vlasov equation to suppress the recurrence.  The
hyperbolic moment method \cite{Fan_new} is first proposed for the
Boltzmann equation, and has been adapted to a lot of fields, including
Vlasov equation \cite{Wang, VM2015}, Wigner equation \cite{Tiao} and
quantum Boltzmann equation. Precisely speaking, the moment method in
kinetic theory was first proposed by Grad \cite{Grad} in 1949, and was
extended into arbitrary order cases in \cite{Fan_new}. However, the
loss of hyperbolicity \cite{Muller} of Grad's moment system limits its
application for a long time.  A recently proposed globally hyperbolic
regularization in \cite{Fan_new} remedies the drawback and yields the
Hyperbolic Moment Equations(HME).  In this paper, we adopt this
globally hyperbolic regularized moment method to solve the Vlasov
equation, and call it HME for short.  HME can be treated as an
``adaptive'' Hermite spectral method in the velocity direction
\cite{dvm2moments}, in which the expansion center is translated by the
local macroscopic velocity and rescaled by the thermal velocity.  The
special transformation on velocity space enhances the efficiency of
HME to approximate the Vlasov equation \cite{Wang}. The numerical
effectiveness of HME on the Landau damping problem has been
demonstrated in \cite{Wang}.

However, as a Eulerian method for Vlasov equation, HME also suffers
from the recurrence phenomenon. Therefore, we would study the
recurrence phenomenon of HME and show that HME shares the problem as
the classical discrete velocity method (DVM). Then filtering is
utilized to suppress the recurrence phenomenon when simulating the
Vlasov equation. Notice that HME is a set of partial differential
equations obtained by approximation of Vlasov equation, which
satisfies a lot of physical properties of Vlasov equation, such as the
conservation of mass, momentum and energy, and Galilean invariance.
This puts forward lots of constrains on the filter.  Moreover, the
authors \cite{Kanevsky2006Idempotent} pointed out that for a given
space and velocity space discretization, different time step lengths
would yield different numerical results due to the different
application times of the filter. That is to say that the solution of a
direct application of filter is time-step dependent, so the limit
equation is not clear. To avoid it, we propose a quasi time-consistent
filter for HME, which also preserves the conservation of mass,
momentum and energy, Galilean invariance and the convergence to the
Vlasov equation with the increasing of the moment number.  To
understand the principle of the filter, we present two viewpoints ---
collision operator viewpoint and artificial dissipation viewpoint ---
to analyze the filter and the resulting system.  In both viewpoints,
the filtered method can be treated as a solver of the Vlasov equation,
thus the damping slope and frequency in the Landau damping problem
would keep unchanged. Numerical simulations show that the filter can
suppress the recurrence phenomenon and demonstrate that the damping
slope and frequency are unchanged by the filter. The filtered method
is also applied to the nonlinear two-stream instability problem to
show its numerical efficiency.  Good agreement with the reference
shows the effectiveness of the filtered method.

The outline of the paper is as follows. In Section
\ref{sec:HMM-Vlasov}, HME for the Vlasov equation is briefly reviewed
and the recurrence phenomenon for HME is studied. In Section
\ref{sec:FilterHME}, the detailed procedure of the filter is presented
and we also provide two viewpoints to analysis the effect of the
filter for HME. In Section \ref{sec:num}, the numerical simulations
are performed to suppress recurrence with our filter in linear Landau
damping problem, and two-stream instability problem is also studied.
The paper ends with a conclusion in Section \ref{sec:conclusion}.

%%% Local Variables: 
%%% mode: latex
%%% TeX-master: "article"
%%% End: 

% vim: tw=70:spell
\section{Hyperbolic Moment System for Vlasov Equation}
\label{sec:HMM-Vlasov}
For the Vlasov-Poisson equations (VP) \eqref{eq:vlasov}\eqref{eq:force}, we introduce
the Maxwellian distribution $f_{eq}$ as
\begin{equation} \label{eq:maxwellian}
  f_{eq}(t,\bx,\bv) =  \frac{\rho(t,\bx)}{[\sqrt{2\pi} u_{th}(t,\bx)]^D}
  \exp\left(-\frac{|\bv-\bu(t,\bx)|^2}{2u_{th}^2(t,\bx)}\right),
\end{equation}
where the parameters $\rho(t,\bx)$, $\bu(t,\bx)$ and $u_{th}(t,\bx)$
denote the density, the macroscopic velocity and the thermal
velocity, respectively, and they are related to the distribution
function by
\begin{equation}\label{eq:rhout}
  \begin{aligned}
    \rho(t,\bx) &=  \int_{\bbR^D}f(t,\bx,\bv)\dd\bv, \\
    \rho(t,\bx)\bu(t,\bx) &= \int_{\bbR^D}\bv f(t,\bx,\bv)\dd\bv, \\
    \frac{1}{2}\rho(t,\bx)|\bu(t,\bx)|^2 + \frac{D}{2}\rho(t,\bx) u_{th}^2(t,\bx) 
    &=  \frac{1}{2}\int_{\bbR^D}|\bv|^2 f(t,\bx,\bv)\dd\bv. 
  \end{aligned}
\end{equation}

With the periodicity condition in $\bx$, and the charge neutrality
consition, VP satisfy the conservation of mass, momentum and total
energy \cite{Cheng2013StudyOC}.  Precisely speaking, multiplying the
equation \eqref{eq:vlasov} by $1$ and $\bv$, direct integration with
respect to $\bv$ and $\bx$ on $\bbR^D \times \Omega$ yields the
conservation of mass and momentum
\begin{equation}
  \od{~}{t}\int_{\bbR^D\times \Omega} 
  \begin{pmatrix}
    1\\
    \bv
  \end{pmatrix}
  f(t,\bx,\bv)\dd\bx \dd\bv = 0,
  \quad  t \in \bbR^+. 
\end{equation}

Multiplying the equation \eqref{eq:vlasov} by $|\bv|^2$ and
integrating by parts, we get the conservation of the total energy for
the system \eqref{eq:vlasov} and \eqref{eq:force}:
\begin{equation}
  \label{eq:energy_conversvation}
  \od{~}{t}\left(\int_{\bbR^D\times
      \bbR^D}f(t,\bx,\bv)|\bv|^2\dd\bx\dd\bv +
    \int_{\bbR^D}|\bE(t,\bx)|^2\dd\bx \right) = 0, \quad  t
  \in \bbR^+.
\end{equation}

\subsection{Hyperbolic moment equations}\label{sec:HME}
The key idea of Grad's moment method is expanding the distribution
function around the Maxwellian $f_{eq}$ into Hermite series as
follows:
\begin{equation} \label{eq:expansion}
  f(t,\bx, \bv) = \sum_{|\alpha| \leqslant M}f_{\alpha}(t,\bx)
  \mH^{[\bu(t,\bx),u_{th}(t,\bx)]}_{\alpha}(\bv),
\end{equation}
where $f_{\alpha}(t,\bx)$ are the expansion coefficients and
$\alpha\in\bbN^D$ is the $D$-dimensional multi-index, and
$|\alpha|:=\sum_{d=1}^D\alpha_d$. The basis functions
$\mH^{[\bu(t,\bx),u_{th}(t,\bx)]}_{\alpha}(\bv)$ are generalized
weighted Hermite functions, defined by
\begin{equation} \label{eq:basis}
  \mH^{[\bu,u_{th}]}_{\alpha}(\bv)
  = (-1)^{|\alpha|} \pd{^{|\alpha|}}{v^{\alpha_1} \cdots \partial
    v^{\alpha_D}} \omega^{[\bu,u_{th}]}(\bv),
\end{equation}
where the weight function $\omega^{[\bu,u_{th}]}(\bv)$ is a Gaussian
function with the form
\begin{equation} \label{eq:weight}
  \omega^{[\bu,u_{th}]}(\bv)
  = \frac{1}{[\sqrt{2 \pi }u_{th}]^D} \exp \left( 
    -\frac{|\bv-\bu|^2}{2 u_{th}^2}
  \right). 
\end{equation}
Noticing the definition \eqref{eq:rhout} of the macroscopic parameters
$\rho$, $\bu$ and $u_{th}$, we obtain
\begin{equation}\label{eq:constraints}
  f_0 = \rho , \quad  
  f_{e_k} = 0, \quad k =1, \cdots , D, \quad 
  \sum_{d=1}^D f_{2e_d} = 0,
\end{equation} 
where $e_d\in\bbN^D$ is the $d$-th unit multi-index, i.e its $d$-th
entry is $1$ and other entries are all zero. By plugging the expansion
\eqref{eq:expansion} into the Vlasov equation \eqref{eq:vlasov}, and
applying the globally hyperbolic regularization in \cite{Fan_new}, we
obtain the globally hyperbolic moment system for the Vlasov equation
{\footnotesize 
  \begin{equation}\label{eq:ms}
    \begin{aligned}
      \pd{f_{\alpha}}{t} &+ %{
      \sum_{j=1}^D \left( u_{th}^2 \pd{f_{\alpha - e_j}}{x_j} 
        +u_j\pd{f_{\alpha}}{x_j}
        +(1-\delta_{M,|\alpha|})(\alpha_j + 1) \pd{f_{\alpha+e_j}}{x_j}
      \right)\\   %}
      +\sum_{d=1}^D \pd{u_d}{t} f_{\alpha-e_d} &+ %{ 
      \sum_{j,d=1}^D \pd{u_d}{x_j} \left( 
        u_{th}^2 f_{\alpha-e_d-e_j} 
        +u_jf_{\alpha-e_d}
        + (1-\delta_{M,|\alpha|})(\alpha_j + 1) f_{\alpha-e_d+e_j} 
      \right) \\ %} 
      + u_{th}\pd{u_{th}}{t} \sum_{d=1}^D f_{\alpha-2e_d} &+ %{
      \sum_{j,d=1}^D
      u_{th}\pd{u_{th}}{x_j} \left( u_{th}^2 f_{\alpha-2e_d-e_j} 
        +u_jf_{\alpha-2e_d}
        + (1-\delta_{M,|\alpha|})(\alpha_j + 1) f_{\alpha-2e_d+e_j}
      \right) \\ %}
      &=\sum_{d=1}^DF_df_{\alpha-e_d},
      \qquad |\alpha|\leqslant M,
    \end{aligned}
  \end{equation} 
}%
where $\delta$ is the Kronecker's delta, and $f_{\alpha}$ is taken as
zero if any component of $\alpha$ is negative. We refer the readers to
\cite{Wang} for the details of the derivation. Noticing
\eqref{eq:constraints}, we collect all the independent variables of
$f_{\alpha}$, $\bu$ and $u_{th}$ as a vector $\bw$, then the system
\eqref{eq:ms} can be written in a quasi-linear form
\begin{equation}\label{eq:ms_quasi}
  \bD(\bw)\pd{\bw}{t}+\sum_{j=1}^D\bM_j(\bw)\bD(\bw)\pd{\bw}{x_j}=\bg(\bw),
\end{equation}
where $\bD\pd{\bw}{t}$ corresponds to the time derivative in
\eqref{eq:ms} while $\bM_j\bD\pd{\bw}{x_j}$ describes the convection
term on the $x_j$ direction, and $\bg$ denotes the right hand side of
\eqref{eq:ms}. The detailed expression can be found in \cite{Wang},
and the form is put here only for the sake of convenience.  We point
out that the moment system \eqref{eq:ms_quasi} is globally hyperbolic,
and have the following results.
\begin{proposition}
  The moment system \eqref{eq:ms_quasi} is globally hyperbolic.
  Precisely, for any unit vector 
  $\bn \in\bbS^{D-1}$, 
  the matrix
  \begin{equation}
    \sum_{j=1}^D n_j \bM_j(\bw)
  \end{equation}
  is real diagonalizable,
  and its eigenvalues for $D>1$ are given as 
  \begin{equation}
    \label{eq:characteristic_speeds}
    \bu(t, \bx) \cdot \bn + \rC{k}{m} u_{th}(t, \bx),
    \qquad 1 \leqslant k \leqslant m \leqslant M+1,
  \end{equation}
  where $\rC{k}{m}$ is the $k$-th zero of the $m$-order Hermite
  polynomial $\He_{m}(x)$ and its eigenvalues for the case $D=1$ are
  $u+\rC{k}{M+1}u_{th}$.
\end{proposition}
The proposition is a fundamental condition for the well-posedness of
HME. We refer readers to \cite{Fan_new} for the proof of the
proposition.  It was pointed out in \cite{Torrilhon2006} that the
characteristic speeds of the moment system can be viewed as the
discretization points of the distribution function. Eq.
\eqref{eq:characteristic_speeds} indicates that the discretization
points of the distribution function for HME are the rescaled zeros of
the Hermite polynomials $\bu\cdot\bn+\rC{k}{m}u_{th}$. The points vary
at different position $\bx$ and different time $t$, which is similar
to the moving mesh method\cite{adjerid1986moving}. In this sense, the
moment method can be viewed as an ``adaptive'' Hermite collocation
method.

\subsection{Velocity space filamentation and recurrence}
\label{Sec:filamentation}

A well-known property of the Vlasov equation, so-called {\it
  filamentation}, is that an initially smooth distribution function
may become increasingly oscillatory in velocity space, with a smaller
and smaller wave number as time evolves. Such oscillations lead to
Landau damping and other kinetic effects, but also make Vlasov
equation hard to resolve numerically. But when the wave number of the
particles is smaller than the minimum wave number that the numerical
method could resolve, the numerical method fails to capture the
filamentation structure of the distribution function.  They eventually
lead to an apparent numerical instability, and sometimes the initial
condition artificially reappears and creates a spurious increase in
the amplitude of the electric field, a phenomenon known as recurrence.

To illustrate the main numerical difficulty caused by the
filamentation for the discretization of the Vlasov equation, we
consider the reduced, free-streaming problem for $1$D case
\begin{equation}
  \label{eq:free_streaming}
  \pd{f}{t} + v \pd{f}{x} = 0, \quad f(0,x,v) = g(x,v),
\end{equation}
then the solution to \eqref{eq:free_streaming} is
$f(t,x,v)= g(x-vt,v)$. The free streaming particle motion shears the
initial phase space as time increases, and the initial spatial
structures are stretched into fine scale structures in velocity space.
Then the derivative of $f$ with respect to velocity
\begin{equation}
  \pd{f}{v} = \left. \pd{g(\zeta,v)}{v} \right|_{\zeta=x-vt}
  - t \pd{g(x-vt,v)}{x}
\end{equation}
grows unbounded. For example, the initial value
\[
g(x,v)= \left(1+ A\cos(kx) \right) \frac{1}{\sqrt{2 \pi}} \exp
\left(-\frac{v^2}{2} \right)
\]
gives the exact solution of $f$ as
\[
f(t,x,v) = \left( 1 + A\cos \left[ k(x-vt) \right] \right)
\frac{1}{\sqrt{2 \pi}} \exp \left(-\frac{v^2}{2} \right).
\]
It can be noted that the ``wavelength" in velocity space is
$\lambda_{v} = 2 \pi/(kt)$, and the distribution $f$ becomes more and
more oscillatory due to the $kvt$ term. Moreover, the macroscopic
variables decay with time, for example, the number density
\[
n_{ext} = \int_{\bbR} f \dd v
= 1+ A \cos(kx)\exp \left(-k^2 t^2/2 \right),
\]
decays super-exponentially fast with time.

According to Nyquist-Shannon sampling theorem, one needs at least two
grid points per wavelength to reconstruct the solution. Thus,
representing the distribution function on an equidistant grid will be
impossible when $t> \pi /(k \Delta v)$, and it leads to numerical
instability and recurrence phenomenon. Precisely, assume that we
resolve velocity space with an equidistant grid
$v_j = j \Delta v, j = 0, \pm 1, \pm 2, \cdots, \pm M$, where $M$ is a
large integer, then calculate the numerical integral
\begin{equation}\label{eq:n_num_dvm}
  n_{num}^E = \sum \limits_{j=-M}^{M} 
  \left[1+A\cos \left( k(x-j \Delta v t) \right)\right]
  \frac{1}{\sqrt{2\pi}}\exp \left(-j^2 (\Delta v)^2 / 2\right) \Delta v.
\end{equation}
It turns out to be periodic in time with periodicity
$\mathrm{T_{recurrence}} = \frac{2 \pi}{k \Delta v}$.
\begin{figure}[htbp]
  \centering
  \subfigure[Discrete velocity method]{%
      \label{fig:adv_recurrence_a}
      \begin{overpic}[width=0.4\textwidth]{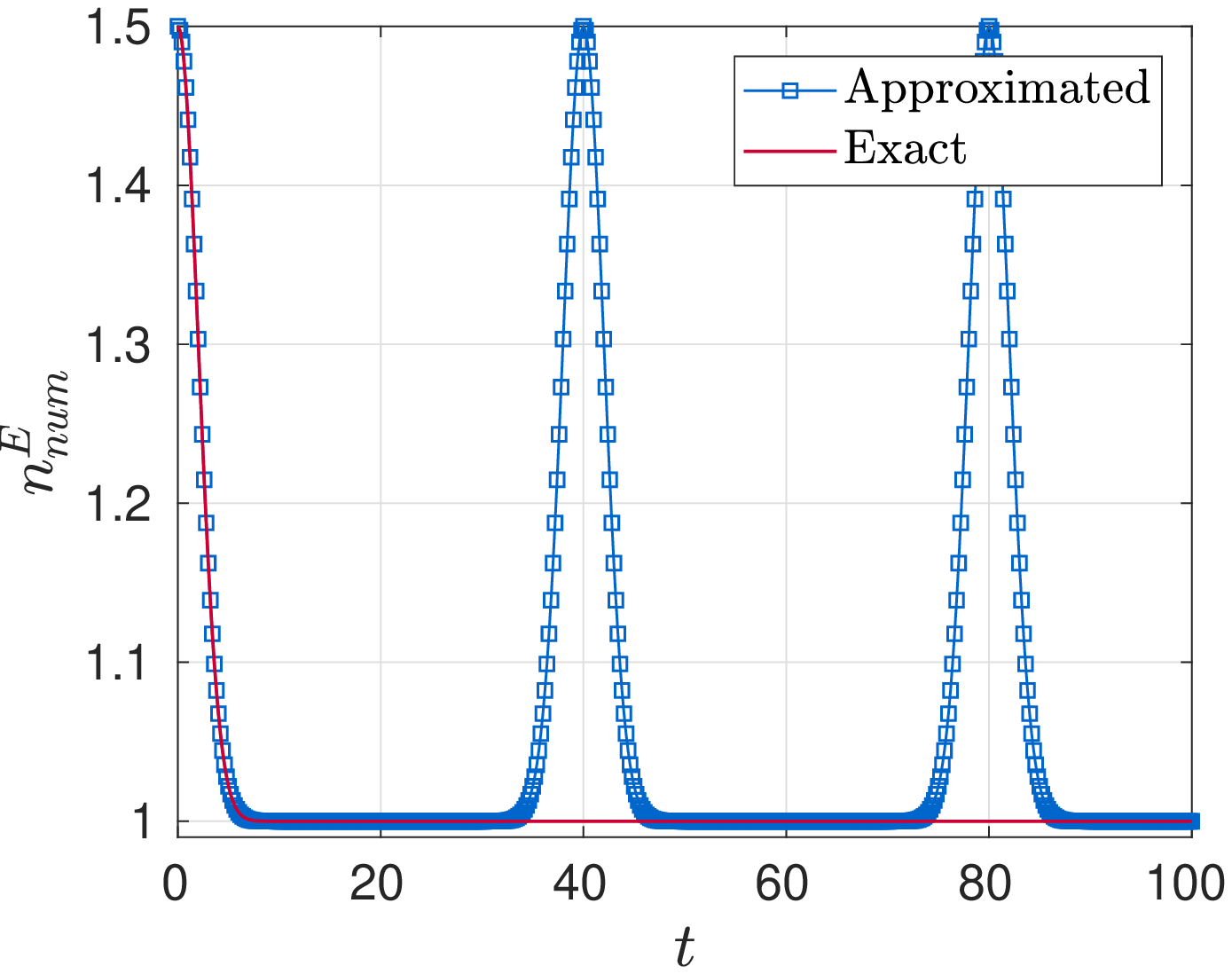}
      \end{overpic}
  }\quad 
  \subfigure[Hermite collocation method]{%
      \label{fig:adv_recurrence_b}
      \begin{overpic}[width=0.4\textwidth]{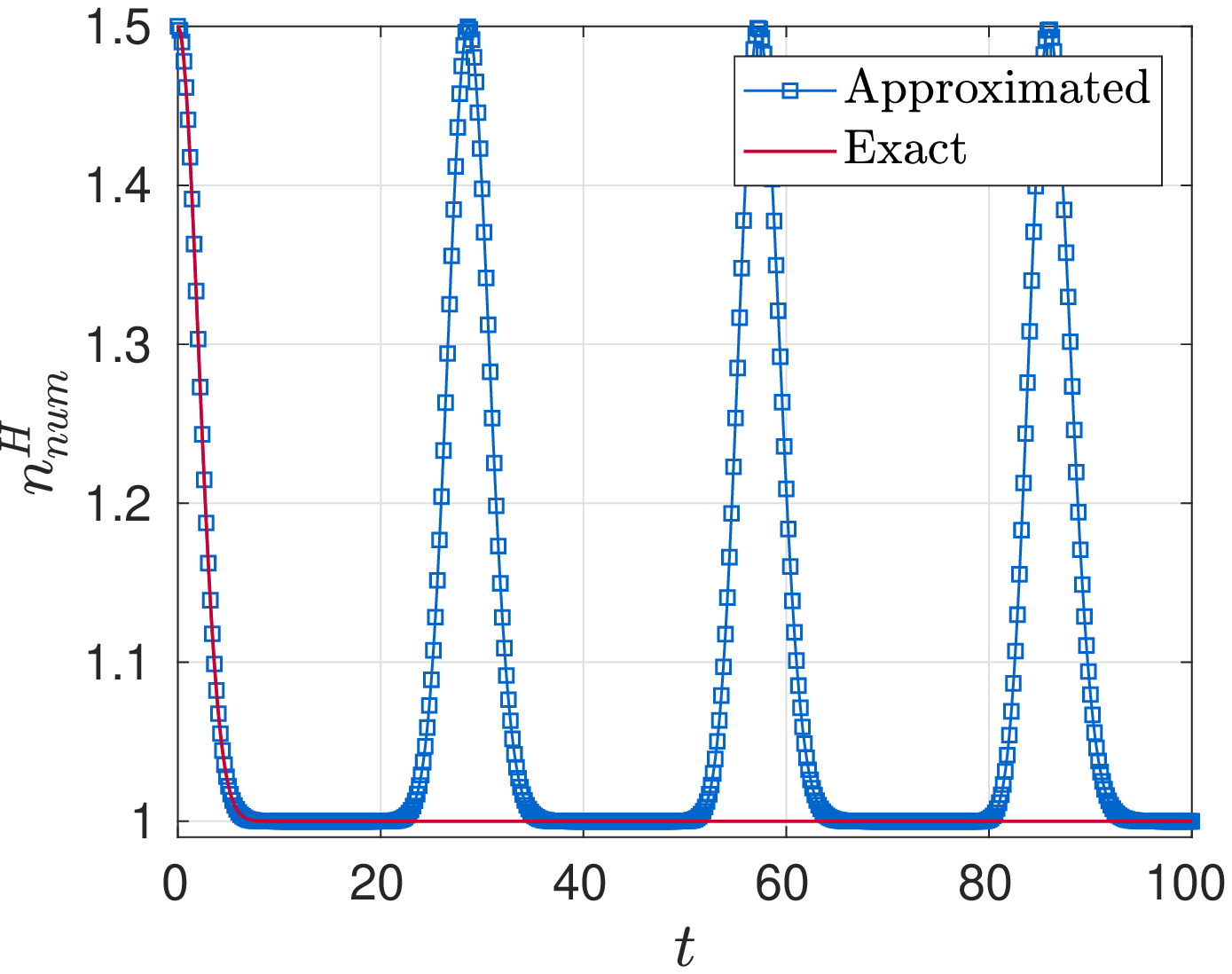}
      \end{overpic}
  }
  \caption{Profiles of the exact and numerical approximations for
    the number density with $A=0.5,k=0.5$. (a) Using the discrete
    velocity method \eqref{eq:n_num_dvm} with $\Delta v=\pi/10$. (b)
    Using the Hermite collocation method \eqref{eq:n_num_H} with
    $M=50$.  }
  \label{fig:adv_recurrence}
\end{figure}
As the recurrence phenomenon shown in Figure
\ref{fig:adv_recurrence_a}, the initial condition reappears
periodically with periodicity $\mathrm{T_{recurrence}} = 40$, while
the exact number density decays super-exponentially. For any
equidistant grid discretization, the recurrence time is proportional
to the grid size for velocity space.

Different from the discrete velocity model, HME expands the
distribution function into the generalized Hermite series
\eqref{eq:expansion} in velocity space, and selects a special set of
characteristic speeds \eqref{eq:characteristic_speeds} such that they
coincide with the Gauss-Hermite interpolation points. As is pointed out
in \cite{Torrilhon2006}, the characteristic speeds can be viewed as a
sort of discretization of the distribution function. Therefore, the
system \eqref{eq:ms} is similar to the discrete velocity model with a
shifted and scaled stencil, thus it also suffers from recurrence effects. It
is not easy to give an accurate estimation  of the recurrence time for
the non-isometric grid discretization. To estimate the recurrence time
for HME, we assume that $u$ is $0$ and $u_{th}$ is $1$. Then the
moment method can be treated as the Hermite collocation method in the
velocity space. Denote the zeros of $\He_{M+1}(v)$ as
$\left\{ v_j \right\}_{j=0}^{M}$, the maximum absolute value of zeros
\[
\max \limits_{j} |v_j|  \sim \sqrt{M},
\]
then the average distance between zeros
\[
\overline{\Delta v} \sim \frac{2\sqrt{M}}{M} = \frac{2}{\sqrt{M}}.
\]
Hence, the recurrence time for the moment method is estimated as
\begin{equation} \label{eq:moment_recurrence_time}
  \mathrm{T_{recurrence}}  
  \sim  \frac{2 \pi}{k \frac{2}{\sqrt{M}}} 
  = \frac{\pi}{k}\sqrt{M}.
\end{equation}
Let $\{w_j\}_{j=0}^M$ be the Gauss-Hermite quadrature weight, then 
\begin{equation}\label{eq:n_num_H}
  n_{num}^{H} =
  \sum_{j=0}^Mw_j[1+A\cos(k(x-v_jt))]\frac{1}{\sqrt{2\pi}}\exp\left(
    -\frac{v_j^2}{2} \right).
\end{equation}
Figure \ref{fig:adv_recurrence_b} shows that the initial condition
also reappears periodically. This result is consistent with the work
in \cite{Landau1946On,UTH, Bourdiec2006Numerical}. It should be noted
that the basis function in \eqref{eq:n_num_H} is Hermite polynomial,
while the basis function in HME is generalized Hermite polynomial.
The upper analysis qualitatively reveals the recurrence phenomenon of
HME and that the recurrence time depends on the square root of the
moment expansion order $M$, which is consistent with the numerical
result in \cite{Wang}.

% vim: tw=70:spell
\section{Filtered Hyperbolic Moment Method}
\label{sec:FilterHME}
As is discussed in the introduction, filtering is a common procedure
to reduce the effects of the Gibbs phenomenon in spectral methods
\cite{canuto2012spectral}. In Section \ref{sec:HME}, we have pointed
out that HME can be treated as an ``adaptive'' Hermite collocation
method for the Vlasov equation, thus it would be natural to apply a
filter on HME to suppress the recurrence. In \cite{Cheng}, it was
shown that the filamentation had little influence on the lower moments
of the distribution function in many cases. Therefore, it is expected
that the filter can suppress the filamentation but only slightly
affect the macroscopic phenomena, for example Landau damping.

However, there are still many unsolved important issues, including 1)
how to choose the filter for HME without destroying its physical
properties, 2) should the filter be applied once, more times or less per
time step, 3) does the filter affect the convergence of HME and 4)
does the filter change the Landau damping rate? In this section, we
will answer these questions one by one.

\subsection{Choosing the filter} \label{sec:filtermethod}

At the first step, we would study the properties of HME and propose
some necessary conditions for the filter to preserve these
properties. Then we are going to give two versions of filters. The
first version is the classical exponential filter
\cite{Hesthaven2008filtering, Gottlieb2001spectral, HouLi2007}.
However, this kind of filter will lead to the potentially
multiplicative net effect of filtering in time
\cite{Kanevsky2006Idempotent}. Noticing this deficiency, a quasi
time-consistent filter, which takes the time step into consideration,
is proposed.

\subsubsection{Necessary conditions of the filter}
\label{sec:condition_filter}
The usual practice of the filter for spectral methods
\cite{parker2015fourier, mcclarren2010robust} is multiplying
spectral coefficients by a filter factor $\sigma$. For HME, the
coefficients $f_{\alpha}$ are corresponding spectral coefficients,
thus the filter can be written as
\begin{equation}
\label{eq:apply_filter}
  \tilde{f}_{\alpha} := f_{\alpha}\sigma\left(\frac{\alpha}{M}\right).
\end{equation}
Accordingly, the distribution function is replaced by 
\begin{equation}\label{eq:filtered_distribution}
  \tilde{f}=\sum_{|\alpha|\leqslant M}\tilde{f}_{\alpha}
  \mH^{[\bu,u_{th}]}_{\alpha}(\bv).
\end{equation}
As is pointed out at the beginning of this section, it is an important
issue to choose the filter without destroying the properties of HME.
Next, we will delineate the necessary conditions for the filter
by studying the properties of HME.

HME is a physical model, derived from the Vlasov equation, thus
it satisfies the Galilean invariance. The filter should
preserve the Galilean invariance of the model, including
the rotational invariant. Hence, we demand the rotational invariant of
the filter, i.e
\begin{equation}\label{eq:filter_rotational}
  \sigma\left(\frac{\alpha}{M}\right)=\sigma\left(\frac{|\alpha|}{M}\right),
  \quad |\alpha|\leqslant M.
\end{equation}
The conservation of mass, momentum and energy are primary properties
of the VP system and HME. Thus the filter is expected to
preserve the conservation laws. Here we demand
\begin{equation}\label{eq:filter_conservation}
  \sigma\left(\frac{|\alpha|}{M}\right)=1,\quad |\alpha|\leqslant M_0,
  \quad \text{and }~ M_0\geqslant2.
\end{equation}
The filter is used to remove the filamentation, which is
caused by the high frequency coefficients. Hence, the filter should be
stronger for higher frequency coefficients, i.e. the filter is a
monotone decreasing function.  Here $\alpha$ is a multi-index, so
there are various definitions of the monotone decreasing. Due to the
rotational invariant relation \eqref{eq:filter_rotational}, we let
\begin{equation}\label{eq:filter_monotone}
  \sigma\left(\frac{|\alpha|}{M}\right) 
  \geqslant \sigma\left(\frac{|\alpha|+1}{M}\right),\quad |\alpha|<M.
\end{equation}
As the order $M$ increasing, the filtered distribution
\eqref{eq:filtered_distribution} should converge to the
distribution, which is a necessary condition for the convergence of the
model. So for any given $\alpha\in\bbN^D$, the strength of the
filter should vanish with the increasing of $M$, i.e.
\begin{equation}\label{eq:filter_limit}
  \lim_{M\rightarrow\infty}\sigma\left(\frac{\alpha}{M}\right)=1.
\end{equation}

All the conditions \eqref{eq:filter_rotational},
\eqref{eq:filter_conservation}, \eqref{eq:filter_monotone} and
\eqref{eq:filter_limit} are necessary conditions for the filter. 
Next, we are trying to explore the filter based on these conditions.

\subsubsection{Exponential filter}
\label{sec:exp_filter}
Exponential filters are widely used in spectral and
pseudo-spectral methods to overcome the Gibbs phenomenon. 
For instance, 
\begin{equation} \label{eq:expoential_filter}
  \sigma(\eta) = \exp\left(-\beta \eta^{\gamma}\right),
\end{equation}
where $\beta = -\ln \epsilon_0$ with $\epsilon_0$ representing the
machine accuracy and $\gamma$ is a constant
\cite{Gottlieb2001spectral}. The stabilizing effect of this filter has
been widely discussed \cite{ Kreiss1979Stability,
  Hesthaven2008filtering, Gottlieb2001spectral}. It has been proved
that this filter is strong enough to stabilize the approximation to
the initial conservation laws and small enough not to ruin the
spectral accuracy of the scheme when the parameters $\beta, \gamma$
are chosen properly. However, the exponential filter fails to
satisfy the condition \eqref{eq:filter_conservation}.

A filter called Hou-Li's filter, is proposed in \cite{HouLi2007} for
Fourier spectral method, which reads
\begin{equation} \label{eq:houli_filter}
  \sigma(\eta) = \left\{
    \begin{array}{ll}
      1,&\text{ if } 0 < \eta \leqslant 2/3,\\
      \exp(-\beta\eta^\gamma),& \text{ if } \eta >2/3,
    \end{array} \right.
\end{equation}
where $\beta=36$ and $\gamma=36$, which subjects $\sigma(1)$ to the
machine precision $\epsilon_0=2^{-53}$. This filter achieved great
success in Fourier spectral method \cite{HouLi2007}, and was also used
in Fourier-Hermite spectral method for Landau damping
\cite{parker2015fourier}. The rotational invariant relation
\eqref{eq:filter_rotational} can be easily satisfied by careful
extension to $D$-dimensional case. A natural extension of
\eqref{eq:houli_filter} is
\begin{equation}\label{eq:filter}
  \sigma\left(\frac{\alpha}{M}\right)=
  \left\{
    \begin{array}{ll}
      1,&\text{ if } |\alpha|/M\leqslant 2/3,\\
      \exp\left(-\beta\left(\frac{|\alpha|}{M}\right)^\gamma\right), 
        & \text{ if } |\alpha|/M>2/3.
    \end{array} \right.
\end{equation}
Moreover, one can easily verify that this filter satisfies the
conditions \eqref{eq:filter_monotone}, \eqref{eq:filter_limit} and
\eqref{eq:filter_conservation}. It is handy to set the parameters
$\beta$ and $\gamma$ here the same as in \eqref{eq:houli_filter}.

\subsubsection{Quasi time-consistent filter}
\label{sec:timeconsitent_filter}
In filtered methods, the filter is usually applied to the distribution
function in each time step. Hence, for a given space and velocity
space discretization, different time step lengths would yield
different numerical results due to the different application times of
the filter. This issue was noticed in \cite{Kanevsky2006Idempotent}
when they compared the result of the explicit Runge-Kutta method and
IMEX-RK method in time discretization.  To avoid this issue, the
authors of \cite{Kanevsky2006Idempotent} proposed a concept of
time-consistent filters. Precisely, they took the time step length
into the filter and denote it as $\sigma(\eta, \Delta t)$.  

\begin{definition}
If a filter $\sigma(\eta, \Delta)$ satisfies 
\begin{equation} \label{eq:time_consistent_filter}
  \sigma(\eta, \Delta t_1)^{k_1} = \sigma(\eta, \Delta t_2)^{k_2},
\end{equation}
when $k_1\Delta t_1 = k_2\Delta t_2$, we would call it a
time-consistent filter. 
\end{definition}

In the following of this section, we will point out that the
time-consistent filter not only fixes the issue of how many times the
filter should be applied on each time step, but also works well in the
numerical simulations of HME for VP.

Here we take the idea of time-consistent into the construction of the
filter and propose a filter as
\begin{equation} \label{eq:filter_t} 
  \sigma\left(\frac{\alpha}{M},\Delta t\right) =
  \left\{ 
    \begin{array}{ll}
      1,&\text{ if } |\alpha|/M\leqslant 2/3,\\
      \exp\left(-\beta \left(\frac{|\alpha|}{M}\right)^\gamma
      g\left(\frac{|\alpha|}{M},\frac{\Delta t}{T_0}\right)\right),
        &\text{ if } |\alpha|/M > 2/3,
    \end{array} \right.
\end{equation}
where
\begin{equation} \label{eq:filter_factor} 
  \beta = 36, \quad \gamma=36,\quad 
  g\left(\eta,\zeta\right) =
  \left(\zeta \right)^{1-\eta^\gamma},
\end{equation}
and $T_0$ is a constant dependent on the dimension of the variables,
and in this paper we always set it as $1$.

\begin{figure}[ht]
  \centering \subfigure[$M=30$]{
      \begin{overpic}[width=.3\textwidth]{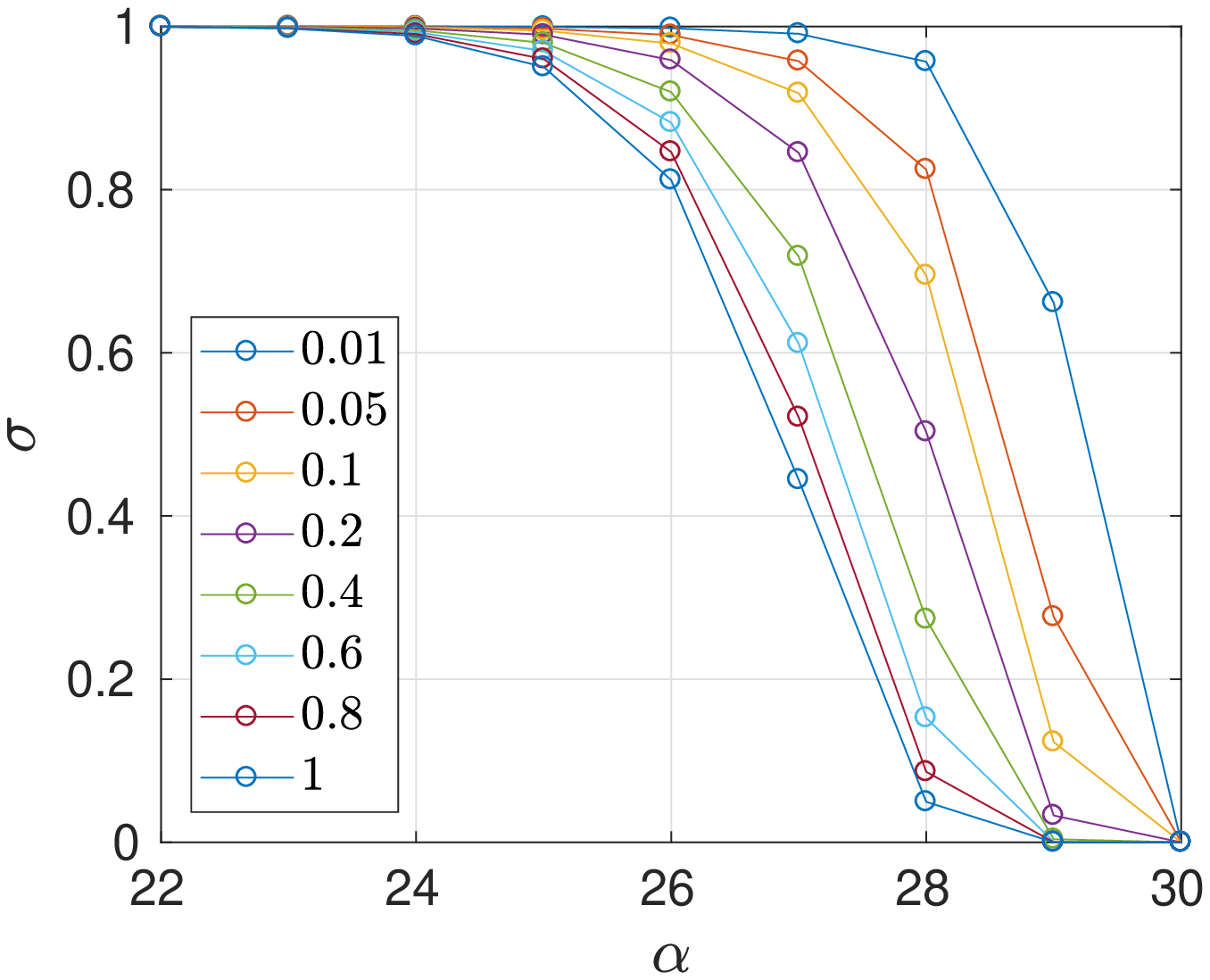}
      \end{overpic}
  } \subfigure[$M=90$]{
      \begin{overpic}[width=.3\textwidth]{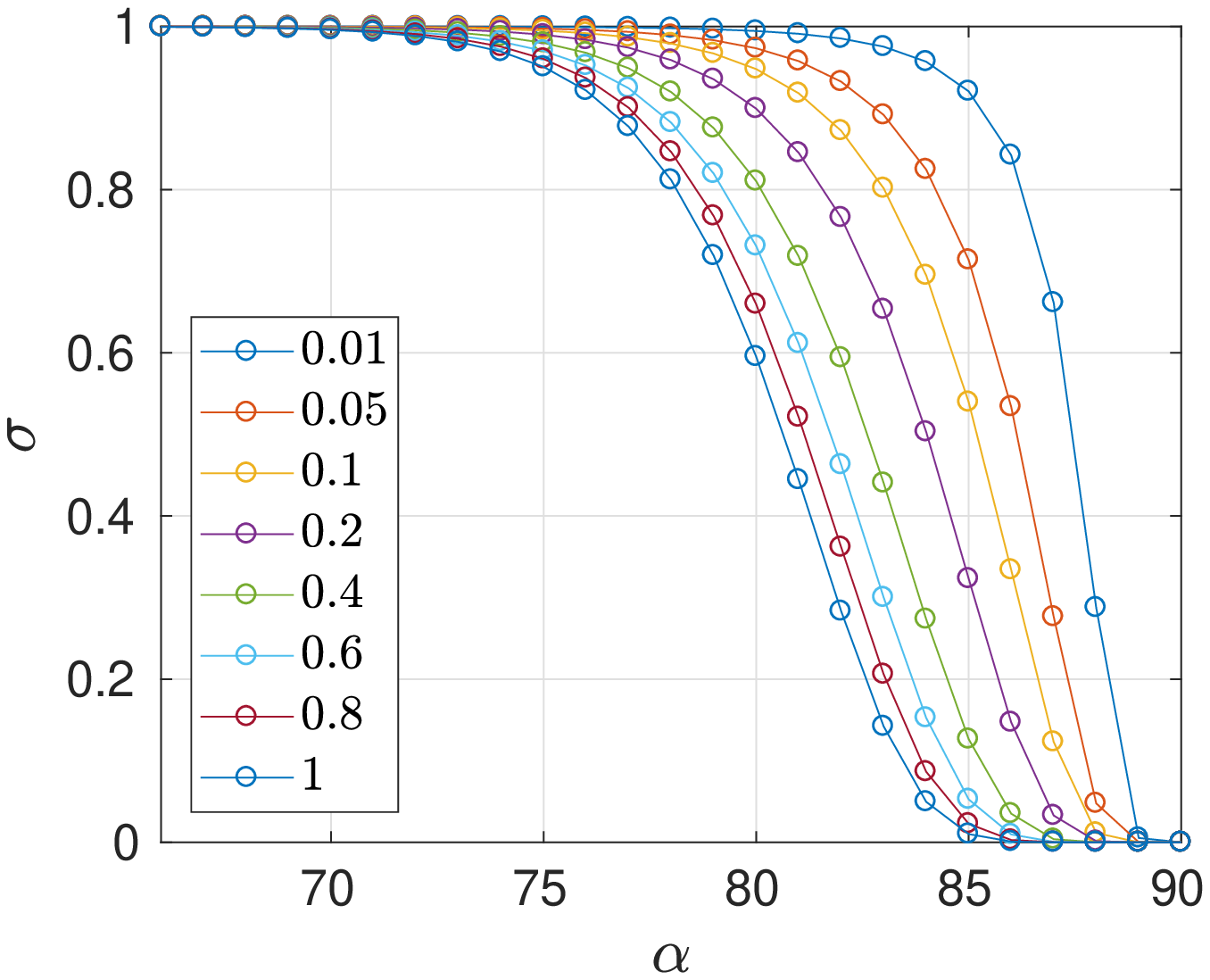}
      \end{overpic}
  } \subfigure[$M=300$]{
      \begin{overpic}[width=.3\textwidth]{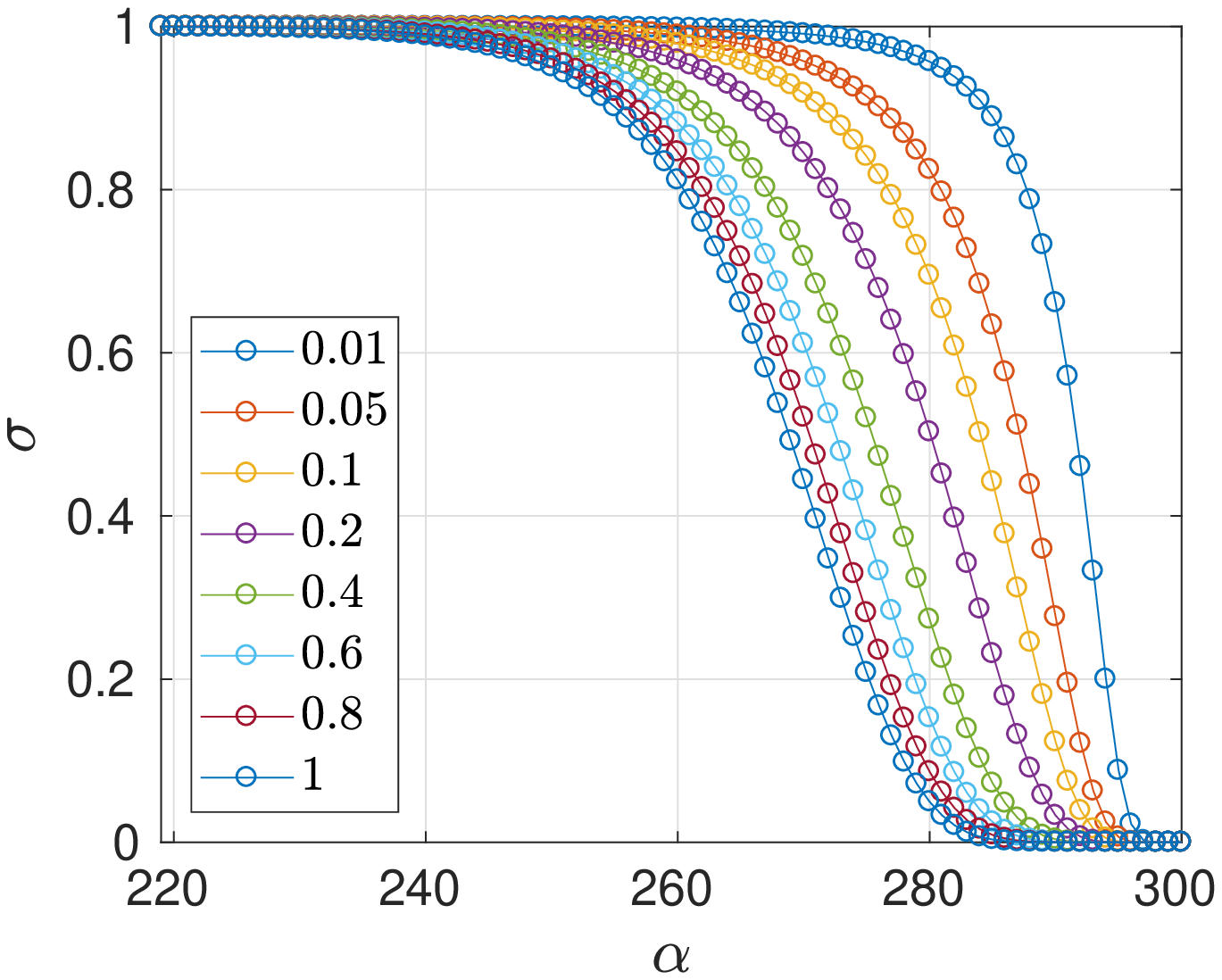}
      \end{overpic}
  }
  \caption{\label{fig:filter}Profiles of the filter
    \eqref{eq:filter_t} with different $\Delta t$ and the order
    $M$. The $x$ and $y$ axis are $|\alpha|$ and
    $\sigma(\alpha/M,\Delta t)$, respectively. The legends are the
    values of $\Delta t/T_0$.}
\end{figure}

If $g(\eta,\Delta t) = 1$, the filter \eqref{eq:filter_t} degenerates
to Hou-Li's filter \eqref{eq:houli_filter}, and if $g(\eta,\Delta
t)=\Delta t$, it is a time-consistent filter. However, for HME, a very
strong filter on the
high-frequency coefficients is necessary to eliminate the
filamentation. To achieve it, one choice is to increase the value
$\beta$, which enhances the strength of the filter for all
$|\alpha|>\frac{2}{3}M$, and another choice is setting
$g(\eta,\Delta t)$ as that in \eqref{eq:filter_factor}, which only
enhances the strength for $|\alpha|\sim M$.  In
Fig. \ref{fig:filter}, we present the filter \eqref{eq:filter_t}
with different $M$ and different relative time step $\Delta t$.

If $|\alpha|\sim M$, $\sigma\left( \frac{\alpha}{M},\Delta t \right)
\sim 0$, so the time consistent condition \eqref{eq:time_consistent_filter} is valid. If
$|\alpha|\sim \frac{2}{3}M$, $g(\eta,\Delta t)\sim \Delta t$, so
the time-consistent conditions is also valid.  Thus we call the
filter \eqref{eq:filter_t} as a quasi time-consistent filter.
Moreover, it is easy to check that the filter \eqref{eq:filter_t}
satisfies all the necessary conditions of the filter
\eqref{eq:filter_rotational}, \eqref{eq:filter_conservation},
\eqref{eq:filter_monotone} and \eqref{eq:filter_limit}.

In this section, we have proposed a quasi time-constant filter to
suppress the recurrence. The time step length will not affect the
total numerical results of the filter. From
Eq. \eqref{eq:apply_filter}, we can see that the filter is applied on
the expansion coefficients of the distribution function instead of
directly on the Vlasov equation. In the following section, we will
show that the effect of the filter on to the Vlasov equation from two
viewpoints, which illustrates that filtering under the framework of
HME is a solver to the Vlasov equation, and filtered HME will predict
the correct physical variables, such as Landau damping rate and
frequency.

\subsection{Mathematical interpretation of the filtering}
In this subsection, we aim to propose two viewpoints to explain the
effect of the filtering on the Vlasov equation.
The first viewpoint explains the filtering as an artificial collision
and the second one interprets the filtering as a dissipative term.
Both viewpoints show that filtered HME is a solver to Vlasov equation,
and the application of filter would not affect the physical variables
people care about.

\subsubsection{Collisional viewpoint} \label{sec:collision} 
In \cite{grant1967fourier}, the authors added an artificial
weakly-collisional operator on Vlasov equation to overcome the
filamentation process. This idea has been well studied in the past
decades \cite{Holloway1996Spectral, ellasson2001outflow,
ng2004complete, camporeale2016velocity}.  Here we would show that the
filter \eqref{eq:filter_t} can also be treated as a kind of artificial
collisional operator.

For simplicity, we take 1D case as an example, and one can extend it
$n$D case without difficulty. By adding the artificial collision operator
\begin{equation} \label{eq:artificialCollision}
  \mC(f) = \nu\left(\sum_{\alpha\leqslant
      M}\nu_{\alpha}f_{\alpha}\mH_\alpha^{[u,u_{th}]}(v) -  f \right)
\end{equation}
on Vlasov equation \eqref{eq:vlasov}, one can obtain the collisional
Vlasov equation
\begin{equation}\label{eq:CollisionValsov}
  \pd{f}{t} + v \pd{f}{x} 
  + F(t,x)\pd{f}{v} = \mC(f).
\end{equation}
Here $\nu_{\alpha}$ and $\nu$ are parameters to be determined, and we
set them as $\nu = \frac{\beta}{T_0}$ and 
\begin{equation}
  \nu_{\alpha} = \left\{ 
    \begin{array}{ll}
      1,  &   \alpha\leqslant M_0,\\
      1-\left(\frac{\alpha}{M}\right)^\gamma,  &   \text{otherwise},
    \end{array} \right.
\end{equation}
in this paper. If let $M_0 = 2$, $\gamma=0$, the operator
\eqref{eq:artificialCollision} degenerates to the classical BGK
collision operator \cite{BGK}. If let $M_0=\lceil \frac{2M}{3}\rceil$,
and we employ the time-splitting scheme to solve the collision term
and other terms, we can obtain the solution of the collision part as
\begin{equation}
  f_{\alpha} \rightarrow \sigma(\alpha/M, \Delta t) f_{\alpha},
  \text{ with }
  \sigma(\alpha/M, \Delta t) = \left\{ 
    \begin{array}{ll}
      1,  &   \alpha\leqslant M_0, \\
      \exp\left( -\beta\left( \frac{\alpha}{M}
      \right)^\gamma\frac{\Delta t}{T_0} \right), &   \text{otherwise},
    \end{array} \right.
\end{equation}
where $\sigma(\alpha/M, \Delta t)$ corresponds to the filter in
\eqref{eq:filter_t} with $g(\eta,\Delta t) = \Delta t$.  From the
viewpoint of ordinary differential equations, we can find a
$\nu_{\alpha}$ which corresponds to the filter in \eqref{eq:filter_t}.
Hence, the filter \eqref{eq:filter_t} can be understood as adding an
artificial collision operator on the Vlasov equation, and the filtered
HME is solving the collisional Vlasov equation
\eqref{eq:CollisionValsov}.

The necessary condition \eqref{eq:filter_limit} for the filter
indicates that for any given $\alpha\in\bbN$, 
\begin{equation}
  \lim_{M\rightarrow\infty}\int_{\bbR}\mC(f)v^\alpha \dd v = 0,
\end{equation}
i.e. $\mC(f)\rightarrow 0$ as $M\rightarrow \infty$, if the Grad's
expansion \eqref{eq:expansion} converges. Particularly, the
collisional Vlasov equation \eqref{eq:CollisionValsov} converges to
the Vlasov equation  as $M\rightarrow \infty$. In other words, the
filtered HME converges to the Vlasov equation as the moment order goes
to infinity. In this sense, the filtered HME is a solver of the Vlasov
equation. Naturally, the filtered HME would predict the correct Landau
damping rate and frequency if $M$ is large enough.

\subsubsection{Dissipative viewpoint}\label{sec:dissipation}
By studying the spectral methods for the hyperbolic systems, the
authors of \cite{Gottlieb2001spectral} pointed out that artificial
dissipation could continuously remove the high-frequency components,
which may help to maintain stability. They also provided a viewpoint
that the filter could also be interpreted as adding the artificial
dissipation. Following the idea in \cite{Gottlieb2001spectral}, we can
obtain the dissipative equation corresponding to the filter
\eqref{eq:filter_t} with $g(\eta,\Delta t) = \Delta t$ for 1D case as
\begin{equation} \label{eq:filtered_vlasov}
  \pd{f}{t} + v \pd{f}{x} + F(t,x) \pd{f}{v} =
  -\beta \frac{(-1)^{\gamma}}{M^{\gamma}} \mathcal{D}^{\gamma} f,
\end{equation}
where the linear differential operator $\mathcal{D}$ is defined by
\begin{equation}
  \mathcal{D}f(t,x,v) = \frac{\partial}{\partial v} \left[
    \exp \left( -\frac{v^2}{2} \right) \frac{\partial}{\partial v}
    \left( \exp \left( \frac{v^2}{2} \right) f(x, v,t) \right)
  \right].
\end{equation}
It is noted that we can also obtain the dissipative equation
corresponding to the filter \eqref{eq:filter_t}, analogously. But the
expression is too complex, so here we take the filter
\eqref{eq:filter_t} with $g(\eta, \Delta t) = \Delta t$ as an example.
It has been proved \cite{Gottlieb2001spectral} that the dissipative
term \eqref{eq:filtered_vlasov} is strong enough to stabilize the
approximation, and yet small enough not to ruin the spectral accuracy
of the scheme for the conservation laws. Further more, it has been
verified that the filter will not affect the damping rate and the
frequency of the Landau damping either. See \cite{Cai2017AFilter} for
more details.

On the other hand, as $M\rightarrow\infty$, the dissipative term
$-\beta\frac{(-1)^\gamma}{M^\gamma}\mathcal{D}^\gamma f\rightarrow0$.
Using the same argument in Section \ref{sec:collision}, it can be
claimed that the filtered HME is a solver of the Vlasov equation and
its solution converges to that of the Vlasov equation as
$M\rightarrow \infty$. Therefore, it would predict the correct Landau
damping rate and frequency if $M$ is large enough.

\begin{remark}
  These two viewpoints have demonstrated that the filtered HME is a
  solver to the Vlasov equation.  They are brought up to show the
  total effect of filtering, other than for the numerical scheme.  The
  numerical algorithm of the filtering is based on
  Eq. \eqref{eq:apply_filter}, which will be described elaborately in
  the next section.
\end{remark}

Before the end of this section, we answer the four questions proposed
at the beginning of this section. In order to choose a filter for HME
without destroying its physical properties, we studied the properties
of HME and proposed four necessary conditions for the filter in
Section \ref{sec:condition_filter}. These four conditions guarantee
the rotational invariant, conservation of mass, momentum and energy
and convergence of HME. The quasi time-consistent filter proposed in
Section \ref{sec:timeconsitent_filter} solves the problem that how
often the filter is applied in each time step. To study the
convergence of the filtered HME and whether the filter changes the
Landau damping rate, we provided two viewpoints: artificial collision
operator in Section \ref{sec:collision} and artificial dissipation in
Section \ref{sec:dissipation} to show that the filtered HME is a
solver of the Vlasov equation and can predict the correct Landau
damping rate and frequency.

%%% Local Variables:
%%% mode: latex
%%% TeX-master: "article"
%%% End:

% vim: tw=70:spell
\section{Numerical Simulations} \label{sec:num} 

In this section, numerical simulations are performed to investigate
the effects of the filtered HME.
We first briefly list the numerical scheme for solving the filtered
HME. Then two classical problems, linear Landau damping and two-stream
instability \cite{Filbet, Wang}, are employed for numerical
simulations. Both problems are set up with periodic boundary
condition and $\rho_0=1$ in \eqref{eq:vlasov}.

\subsection{Numerical scheme}
As is discussed in Section \ref{sec:FilterHME}, the filter is applied
once in each time step. In this subsection, we will first briefly
introduce the numerical scheme to solve VP and then list the outline
of the whole numerical scheme. 

\subsubsection{Numerical scheme to solve VP}
To solve VP, the numerical scheme in \cite{Wang} is employed with
filters added in. We refer readers to \cite[Section 3]{Wang} for more
details of the numerical scheme, and only a brief description of the
numerical scheme is listed here. By a standard fraction step method,
we split VP into the convection step and the acceleration step. From
the deduction in \cite{Wang}, we can find that the acceleration step
only contains the electric force $\bF$ in the governing
equations. Thus VP is split as
\begin{itemize}
\item the convection part:
\begin{equation}
\label{eq:convection_part}
\pd{f}{t} + \bv \cdot \nabla_{\bx}f =  0,
\end{equation}
\item the acceleration part: 
  \begin{equation}
  \label{eq:force_part}     
  \begin{gathered}
    \partial_t\bu = \bF,\\
    \bF(t,\bx,\bv) = \boldsymbol{E}(t,\bx), \quad \boldsymbol{E}(t,\bx)
    = -\nabla_{\bx} ~\phi(t,\bx),\quad -\Delta_{\bx} \phi = \rho(t, \bx)
    - \rho_0.
  \end{gathered}
\end{equation}
\end{itemize}

Here we restrict our study in the 1D spatial space.  The standard
finite volume discretization is adopted in the $x$-direction.  Suppose
$\Gamma_h$ to be a uniform mesh in $\bbR$, and each cell is identified
by an index $j$. For a fixed $x_0 \in \bbR$ and $\Delta x > 0$,
\begin{equation}
  \Gamma_h = \big\{T_{j} = x_0 + \left( j\Delta x, ~(j + 1)\Delta x
  \right): j \in  \mathbb{Z}\big\}.
\end{equation} 
The numerical solution which is  the approximation of the distribution
function $f$ at  $t = t_n$ is denoted as
\begin{equation}
  f_h^n(x, \bv) = f_j^n(\bv)=\sum_{|\alpha| \leqslant M} f_{j, \alpha}^n
\mathcal{H}_{\alpha}^{[\bu_j^n, {u_{th,j}^n}]}(\bv),  \quad x \in T_j.
\end{equation}

For the convection part, the distribution function is updated by the
conservative part and the regularization part (see
\cite[Eq. (3.9)]{Wang})
\begin{equation}
\label{eq:scheme}
  f_{j}^{n+1,\ast}(\bv) = f_{j}^n(\bv) + K_{1,j}^n(\bv) + K_{2,j}^{n}(\bv).
\end{equation}
Here, $K_{1,j}^n$ is discretized in the conservative formation as
\begin{equation}
\label{eq:con}
K_{1,j}^n(\bv) = -\frac{\Delta t^n}{\Delta x}
\left[F_{j+\frac{1}{2}}^n(\bv) - F_{j-\frac{1}{2}}^n(\bv)\right],
\end{equation}
where $F_{j+\frac{1}{2}}^n$ is the numerical flux between cell $T_{j}$
and $T_{j+1}$ at $t^n$ and the same HLL scheme \cite[Eq. (3.11)]{Wang}
is utilized here. Similarly, the numerical approximation for the
regularization part $K_{2,j}^n$ is
\begin{equation}
\begin{split}
\label{revision_scheme}
K_{2,j}^n(\bv) & = -\frac{\Delta t}{2\Delta x}\sum_{|\alpha|=M}
\left(\alpha_1+1\right)\sum\limits_{d=1}^3\left(f_{\alpha-e_d+e_1}^{n}
   \left(u_{d,j +1}^n - u_{d,j-1}^n\right)  \right. \\
 & \qquad \qquad  \left. + f_{\alpha-2e_d+e_1}^{n}u_{th, j}^n\left(u_{th,j + 1}^n -
     u_{th,j-1}^n\right) \right)
 \mathcal{H}_{\alpha}^{[\bu_j^n, u_{th,j}^n]}
\left(\bv\right), \quad |\alpha| = M.
\end{split}
\end{equation}

In the 1D spatial space case, the acceleration part is approximated as
\begin{equation} \label{eq:force_solve}
   u^{n+1}_{1,j} = u^{n+1,\ast}_{1,j} 
  + \Delta t F_{1,j}^{n+1},
\end{equation}
where $u_{1,j}^{n+1,\ast}$ and $F_{1,j}^{n+1}$ is the first entry of
the macroscopic velocity $\bu$ and the electric force $\bF$ in the
$j$-th cell after the convection step at $t = t^n$, respectively.  The
electric force $E_{1, j}^{n}$ is updated as
\begin{equation}
  \label{eq:E}
  \begin{aligned}
    - \frac{\phi_{j+1}^{n+1} - 2\phi_j^{n+1} +
      \phi_{j-1}^{n+1}}{\Delta x^2}
    =  \rho_j^{n+1} - \rho_0 , \qquad 
    F_{1,j}^{n+1} =  E_{1,j}^{n+1} = - \frac{\phi_{j+1}^{n+1} -
      \phi_{j-1}^{n+1}}{2\Delta x}.
  \end{aligned}
\end{equation}
where $\rho_j^{n+1}$ is the density in the $j$-th cell after the
convection step, for the reason that the density is not updated in the
acceleration step and the collision step. 

By now, we have introduced the algorithm to solve VP. The application
of filter can be treated as the revision of the moments, and the
detailed application of the filter is explained in the outline of the
algorithm in the following.

\subsubsection{Outline of the algorithm}
The outline of the algorithm is as follows:
\begin{enumerate}
\item Let $n=0$, $t=0$. Set the moment order $M$ and the initial value
  $f_{j,\alpha}^n$, $\bu_j^n$ and $u_{th,j}^n$ on the $j$-th mesh
  cell;
\item Calculate the time step length $\Delta t^n$ according to the
  CFL condition
  \begin{equation}\label{eq:CFL}
      \Delta t^n = \mathrm{CFL}\frac{\Delta x}{\lambda_{\mathrm{max}}},
      \quad 
      \lambda_{\mathrm{max}}=\max_{j}\left(|\bu_j^n|+\rC{0}{M+1}u_{th,j}^n\right),
  \end{equation}
  where $\mathrm{CFL}$ is the CFL number and $\rC{0}{M+1}$ is the
  maximum zeros of $\He_{M+1}(v)$;
\item Use the numerical scheme in the last section to solve the Vlasov
  equation \eqref{eq:vlasov} by one time step, and denote the solution
  in $i$-th mesh cell as $f_{j,\alpha}^{n,*}$, $\bu_j^{n,*}$ and
  $u_{th,j}^{n,*}$;
\item Filtering:
  $f_{j,\alpha}^{n+1}\leftarrow \sigma\left(\frac{\alpha}{M},\Delta
    t^n\right)f_{j,\alpha}^{n,*}$,
  $\bu_j^{n+1}\leftarrow \bu_j^{n,*}$,
  $u_{th,j}^{n+1}\leftarrow u_{th,j}^{n,*}$;
\item Let $n \leftarrow n+1$, $t\leftarrow t+\Delta t^n$ till
  $t\geqslant t_{end}$.
\end{enumerate}
In this section, we always set $\mathit{CFL}=0.45$ if not stated
specifically.

\subsection{Linear Landau damping}
For VP, one important property people care about is the time evolution
of the square root of the electric energy, which is defined as
\begin{equation} \label{eq:energy} 
  \mathcal{E}(t) = \left(\int_{\Omega}|\bE(t, \bx)|^2\dd \bx\right)^{1/2}.
\end{equation}
According to Landau's theory, the time evolution of $\mathcal{E}(t)$
is expected to be exponentially decaying with a fixed rate $\gamma$
and fixed frequency $\omega_R$, which are given by the dispersion
relation.  Another property people care about is the conservation of
mass, momentum and energy. As is proved in \cite{Wang}, this numerical
scheme keeps the conservation of total mass and momentum. The total
energy is defined as
% \begin{equation}
%   \label{eq:total_mom}
%   \begin{split}
%     &\mathcal{M}_{total}(t) = \int_{\Omega} \rho(t, \bx)\dd \bx, \\
%     &\mathcal{U}_{total}(t) = \int_{\Omega} \rho(t,\bx)\bu(t,\bx) +
%     \bE(t,\bx)\dd \bx.
%   \end{split}
% \end{equation}
% and the total energy, which includes the electric energy, the kinetic
% and internal energy of the particles, is defined as 
\begin{equation}
  \label{eq:total_energy}
  \mathcal{E}_{total}(t) = \int_{\Omega} |\bE(t, \bx)|^2 + \rho(t,\bx) |\bu(t,\bx)|^2
  + \rho(t,\bx) u_{th}^2(t, \bx) \dd \bx.
\end{equation}
The variation of the total energy will be studied in the numerical
simulations.  In the following, we first study the linear Landau
damping problem for 1D case to demonstrate the effectiveness of the
filter and then study the 2D case to show its validity for
high-dimensional case.

\subsubsection{1D linear Landau damping} \label{sec:1d} 
In this subsection, we study the filtered HME by 1D linear Landau damping.
Same initial state as that in \cite{Filbet,Wang} is adopted 
\begin{equation} \label{eq:initial_data}
  f(0,x,v) = \frac{1}{\sqrt{2\pi}}e^{-v^2/2}(1+A \cos(kx)), \quad
  (x,v) \in (0,L)\times \bbR,
\end{equation}
where $A = 10^{-3}$ is the amplitude of the perturbation, $k$ denotes
the wave number, and the periodic length is $L = 2\pi/k$.

Next, we study the properties of the filtered HME in detail.

\begin{figure}[!ht]
  \centering
  \psfrag{E}{\footnotesize $\mathcal{E}$}
  \subfigure[$M=50$]{
    \includegraphics[width=0.4\textwidth]{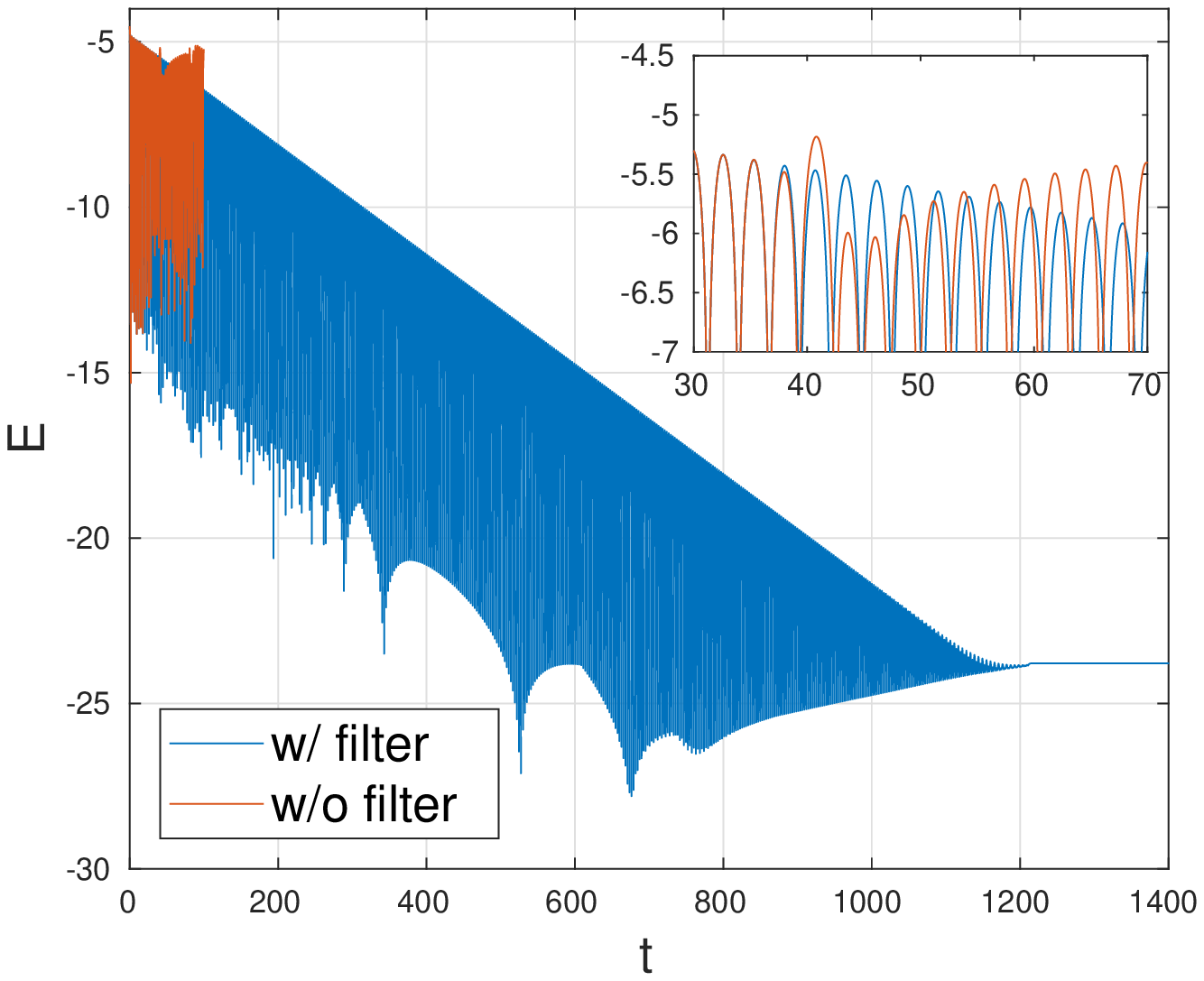}
  }\quad
  \subfigure[$M=80$]{
    \includegraphics[width=0.4\textwidth]{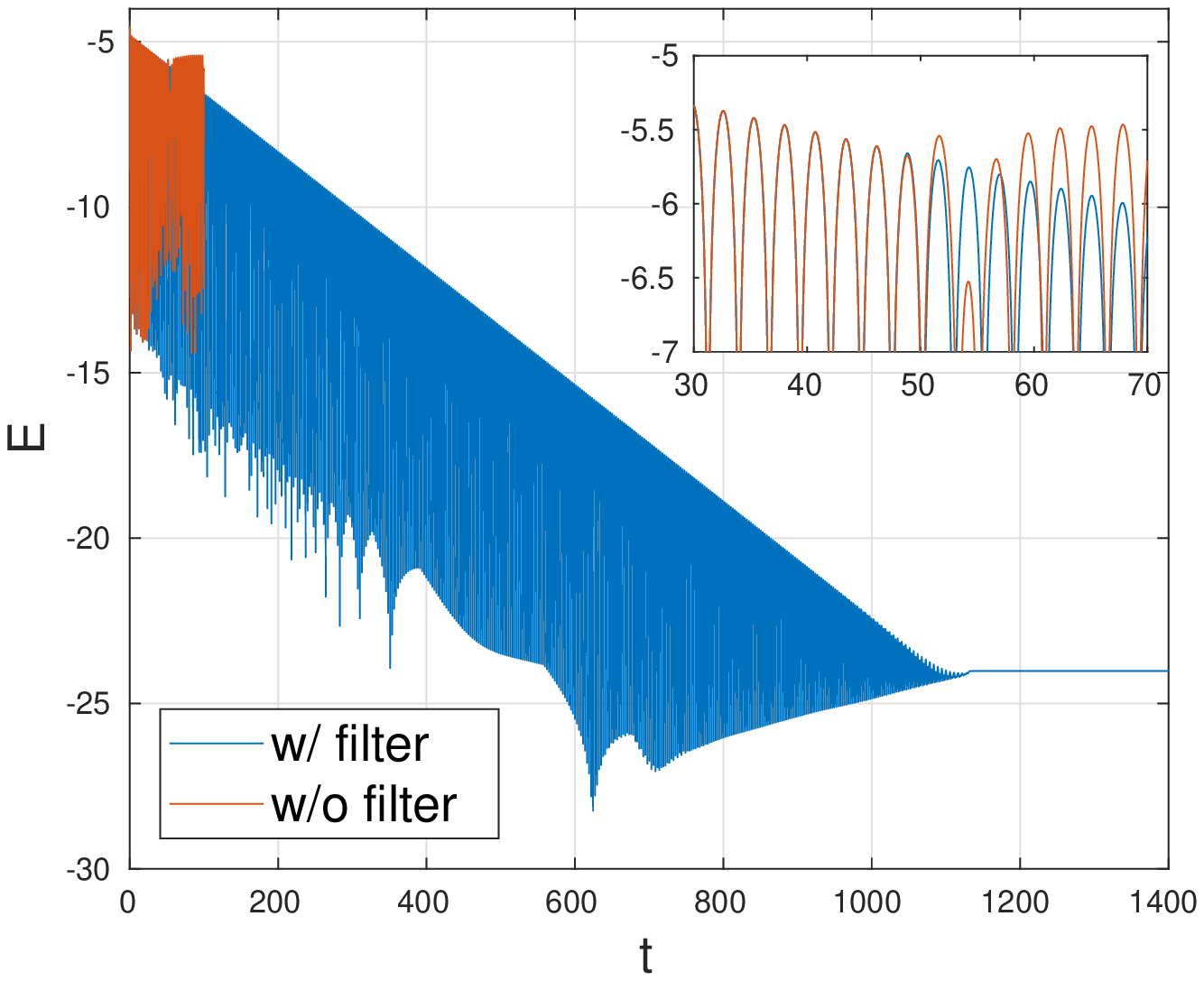}
  }
  \caption{ \label{fig:spatial_03_total}
  Time evolution of $\ln(\mE(t))$ with $N=3200$, $k=0.3$ for different
  $M$. The blue line is that with filter while the red line is that
  without filter.}
  %Exponential damping of $\mathcal{E}(t)$ on the spatial
  %  grids $N = 3200$ and the number of moments $M=50, 80$
  %  respectively, with the wave number $k = 0.3$. The curves are
  %  profile of $\ln(\mathcal{E}(t))$. The blue line is that with
  %  filter while the red line is that without filter.}
\end{figure}

\paragraph{Long time behavior}
The recurrence phenomenon of HME breaks its Landau damping. For the
filtered HME, it is expected that the filter can restrain the
recurrence phenomenon, so the Landau damping can sustain for a long
time. By performing simulations for both HME and the filtered HME with
the wave number $k = 0.3$, the grid number $N=3200$ and the number of
moments $M = 50$ and $80$, we present the time evolution of $\mE(t)$
in Figure \ref{fig:spatial_03_total}.  One can observe that for HME,
the recurrence appears in a short time, but for the filtered HME, the
damping of the energy sustains for a long time until the solution
reaches the machine precision.  Hence, the filter HME has a good
behavior on restraining the recurrence.

\paragraph{Quasi time-consistent filter} 
\begin{figure}[!ht]
  \centering
  \psfrag{E}{\footnotesize $\mathcal{E}$}
  \begin{overpic}[scale=.5]{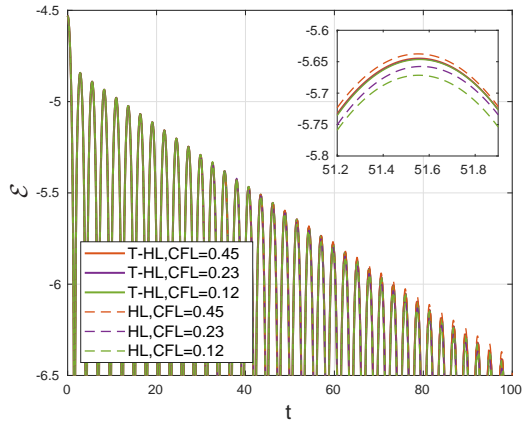}
  \end{overpic}
  \caption{ \label{fig:CFL_03}
  Time evolution of $\ln(\mE(t))$ for the filtered method with
  Hou-Li's filter (HL) and quasi time-consistent filter (T-HL)
  \eqref{eq:filter_t} with different CFL numbers.}
  %Exponential damping of $\ln(\mathcal{E}(t))$ 
  %  for different CFL numbers with the wave number $k=0.3$. ``HL'' and
  %  ``T-HL'' stand for the Hou-Li's filter \eqref{eq:houli_filter} and
  %  quasi time-consistent filter \eqref{eq:filter_t}, respectively.}
\end{figure}
The quasi time-consistent filter is proposed to reduce the effect of
the filter with respect to the time step. For different time step with
the spatial discretization unchanged, which corresponds to different
CFL numbers, the numerical results of the filtered HME should be
almost same. But for Hou-Li' filter, which is employed to solve VP in
\cite{parker2015fourier}, if the CFL number is halved, the filter is
applied twice in the previous one time step. Therefore, its numerical
results may be quite different. Figure \ref{fig:CFL_03} presents the
time evolution of the energy $\mathcal{E}$ for three different CFL
numbers as $\rm{CFL}_1 = 0.45$, $\rm{CFL}_2 = 0.23$, and
$\rm{CFL}_3 = 0.12$.  The setup is $N=3200$, $k=0.3$ and $M=50$. The
numerical results support our conjecture. Precisely speaking,
different CFL numbers barely change the behavior of $\mathcal{E}$ for
the filtered method with filter \eqref{eq:filter}, while the method
with Hou-Li's filter gives different damping effects of
$\mathcal{E}$ under different CFL numbers.

\begin{figure}[!ht]
  \centering
  \psfrag{E}{\footnotesize $\mathcal{E}$}
  \subfigure[Landau damping phenomenon]{
    \includegraphics[width=0.45\textwidth]{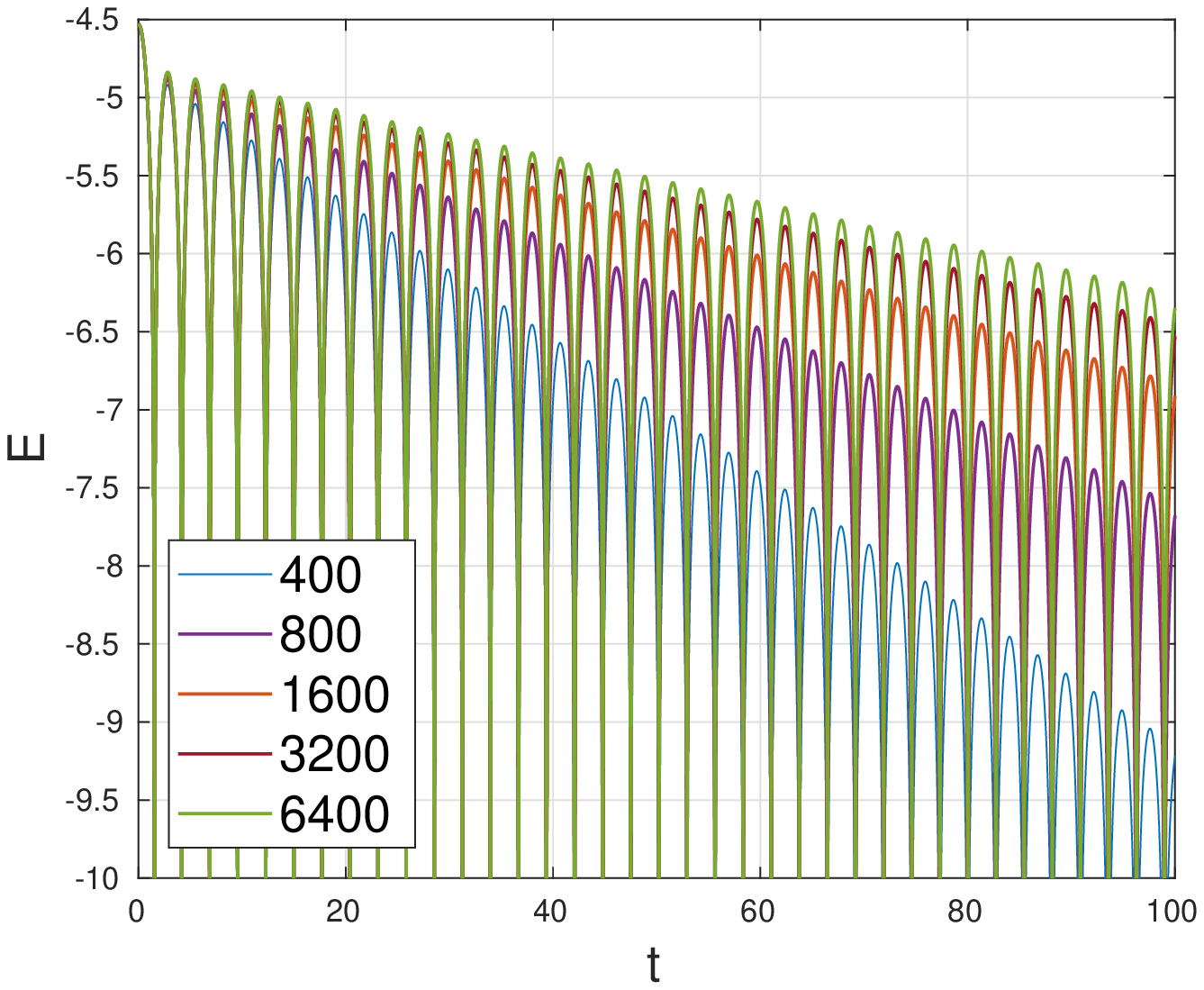}
  }\quad
  \subfigure[Landau damping rate]{
    \includegraphics[width=0.45\textwidth]{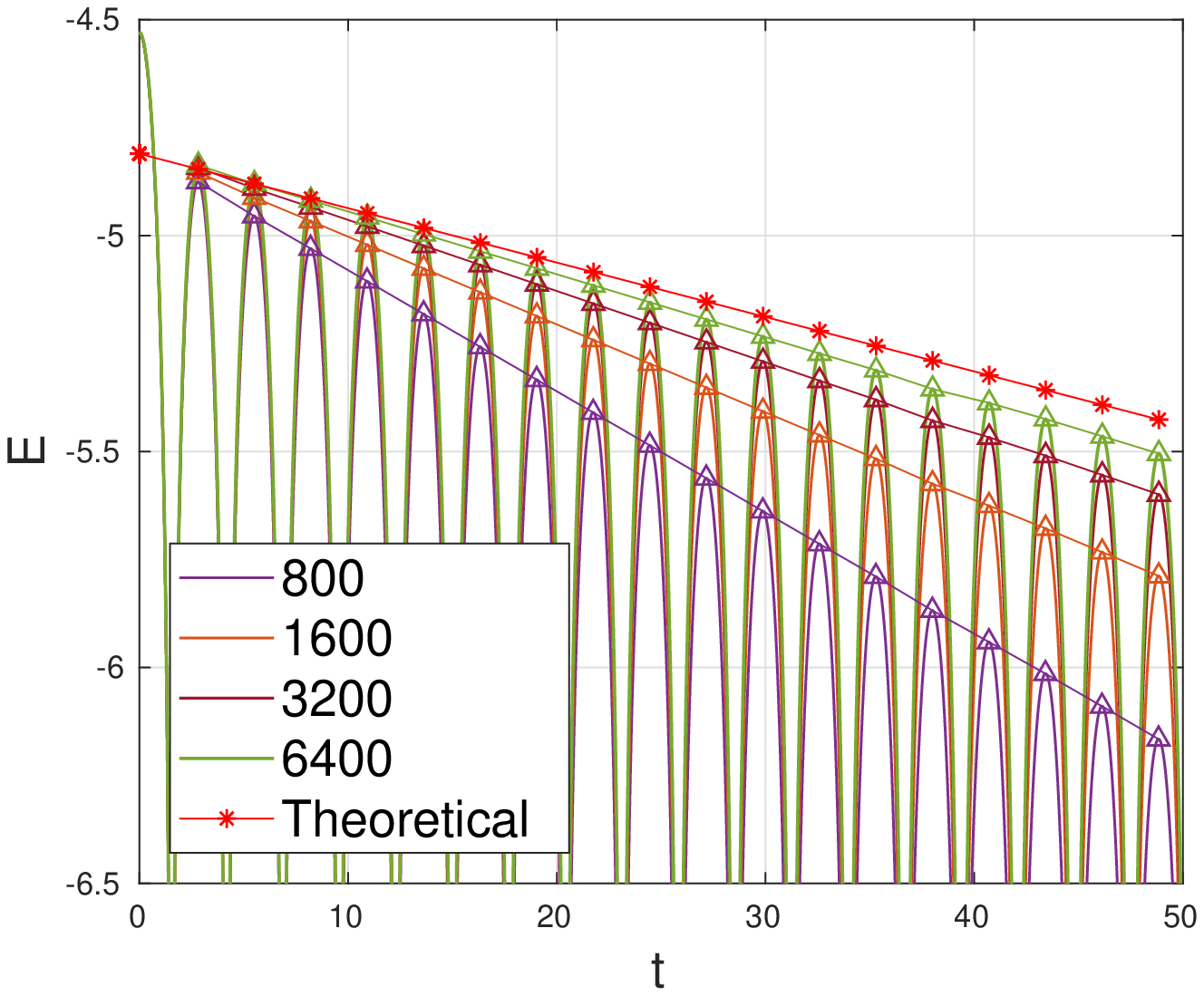}
  }
  \caption{\label{fig:spatial_03} 
  Time evolution of $\ln(\mE(t))$ for different spatial grid steps
  with $k=0.3$ and $M=50$. The slopes of the curves are the numerical
  damping rates calculated by the least square fitting of the peak
  value points of $\mathcal{E}(t)$.  
  }
    %Exponential damping of $\mathcal{E}(t)$ on different spatial grids
    %with $k = 0.3$ and $M=50$. The curves are $\ln(\mathcal{E}(t))$ on
    %different spatial grid sizes. The slopes of the curves are the
    %numerical damping rates by the least square fitting of the peak
    %value points of $\mathcal{E}(t)$. The slope of the star line is
    %the damping rate given by the theoretic data.
\end{figure}

\begin{figure}[!ht]
  \centering
  \subfigure[Landau damping rate $\gamma$]{
    \label{fig:slope_frequency_03_rate}
    \includegraphics[width=0.45\textwidth]{slope_k_03M_50_c.eps}
  }\quad 
  \subfigure[Landau damping frequency $\omega_R$]{
    \label{fig:slope_frequency_03_fre}
    \includegraphics[width=0.45\textwidth]{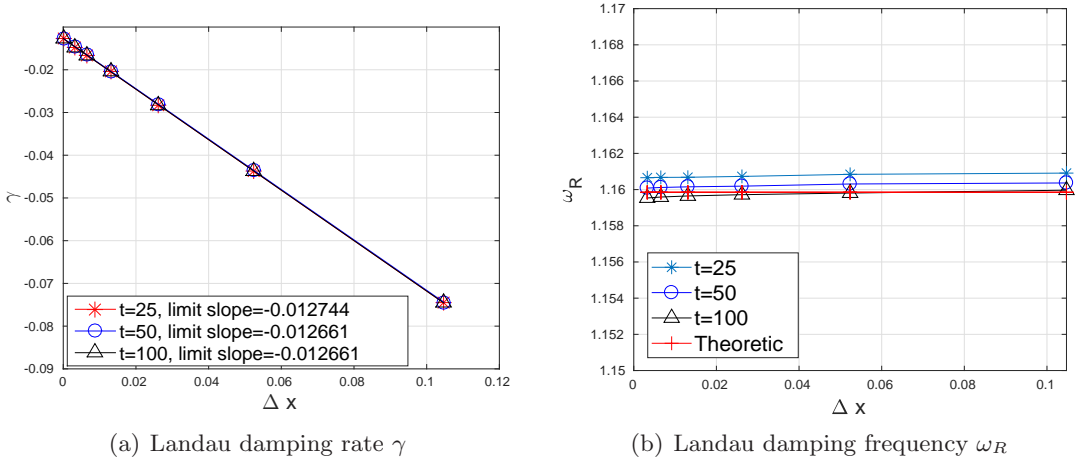}
  }
  \caption{\label{fig:slope_frequency_03} Profile of the numerical
    damping rates and the numerical frequencies with respect to the
    space grid step $\Delta x$ with $k=0.3$ and $M=50$.  The total
    evolution time is $t = 25, 50$ and $100$ respectively.  The
    damping rate $\gamma$ in (a) is obtained by the least square
    fitting of the damping rate with different $\Delta x$. The
    intercept of the line on $y$-axis is the limit damping rate.  The
    frequency $\omega_R$ of the electric field in (b) is estimated by
    counting the peaks in Figure \ref{fig:spatial_03}.}
%  \caption{(a) The linear dependence of the numerical damping
%    rates in spatial grid size with the wave number $k=0.3$. The
%    evolution time is $t = 25, 50$ and $100$ respectively. The $x$-axis
%    is the spatial grid size $\Delta x$ and the $y$-axis is numerical
%    damping rate. The line is obtained by the least square fitting of
%    the damping rate with different $\Delta x$. The intercept of the
%    line on $y$-axis is the limit damping rate. (b) The numerical
%    frequencies with different spatial grid sizes with wave number
%    $k=0.3$. The evolution time is $t = 25, 50$ and $100$
%    respectively. The $x$-axis is the spatial grid size $\Delta x$ and
%    the $y$-axis is numerical frequency. The frequency of the electric
%    field is estimated by counting the peaks in Figure
%    \ref{fig:spatial_03}.}

\end{figure}

\paragraph{Damping rate and frequency}
As argued in Section \ref{sec:collision} and \ref{sec:dissipation},
the filtered HME is a solver of VP, so the filtered HME
would predict the correct Landau damping rate and frequencies if the
moment order $M$ is large enough. To study it, we first study the
spatial convergence to avoid the effects of the spatial error to the
moment convergence and then study the moment convergence.
\subparagraph{Spatial convergence}
\label{para:spatial}
Figure \ref{fig:spatial_03} presents detailed illustrations of the
evolution of $\mathcal{E}(t)$ with different grid sizes as $N = 200$,
$400$, $800$, $1600$, and $3200$.  As the grid number increasing, the
damping rate converges.  The least square fitting method used in
\cite{Wang} is utilized here to approximate the final damping
rate. The analytical solution of Landau damping rate and frequency for
$k=0.3$ is
\begin{equation}
  \label{eq:dispersion_03}
  \omega_R = 1.1598, \quad \gamma = -0.0126.
\end{equation}

Figure \ref{fig:slope_frequency_03_rate} presents the damping rates
for different spatial grid steps and different total evolution time.
The numerical damping rates are in a linearly and monotonically
converging pattern with the spatial grid size $\Delta x$ going to zero
and the limit slope we obtained is in perfect agreement with the
theoretical data. The Landau damping frequencies in Figure
\ref{fig:slope_frequency_03_fre} also has a good convergence to the
theoretical results with the increasing of spatial grids.

Moreover, the damping rates and frequencies for different end time
$t_{end}$ are also studied. Three end time $t_{end}=25$, $50$ and
$100$ correspond to the cases before the recurrence, just after the
recurrence and a long time after the recurrence of HME. Figure
\ref{fig:slope_frequency_03_rate} shows that the damping rates before
and after the recurrence of the HME are almost same, while Figure
\ref{fig:slope_frequency_03_fre} shows that the frequencies before and
after the recurrence of the HME are just slightly different. This
indicates that the filtered HME can predict the correct damping rate
and frequency for a long time after the recurrence of HME.

\begin{figure}[!ht]
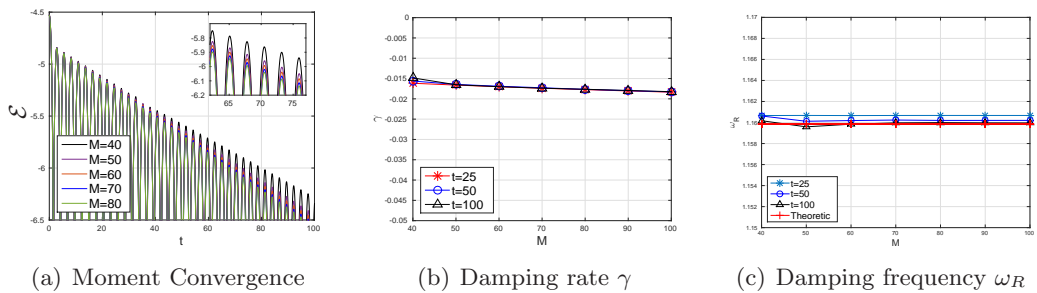

  \centering
  \psfrag{E}{\footnotesize $\mathcal{E}$}
  \subfigure[Moment Convergence]{
    \label{fig:moment_03_a}
    \includegraphics[width=0.27\textwidth]{k_03N_3200M_50_moment_c.eps}
  }\quad 
  \subfigure[Damping rate $\gamma$]{
    \label{fig:slope_frequency_03_moment_rate}
    \includegraphics[width=0.27\textwidth]{slope_k_03N_3200_c.eps}
  }\quad 
  \subfigure[Damping frequency $\omega_R$]{
    \label{fig:slope_frequency_03_moment_fre}
    \includegraphics[width=0.27\textwidth]{k_03N_3200_aveperiod_c.eps}
  }
  \caption{\label{fig:moment_03}
  (a) Time evolution of $\ln(\mE(t))$ for different $M$ with 
  $N=3200$, $k=0.3$. (b), (c) Numerical Landau damping rates 
  and numerical Landau damping frequencies for different 
  $M$ with $k=0.3$, $N=3200$ and  different total evolution time, 
  respectively.
  }
  %(a) Exponential damping of $\mathcal{E}(t)$ on different
  %  number of moments $M$ with the wave number $k=0.3$ and the gird
  %  size $N=3200$. The curves are $\mathcal{E}(t)$ in logarithm scale
  %  form for different number of moments.  (b) The linear dependence
  %  of the numerical damping rates with different number of moments
  %  with the wave number $k=0.3$. The evolution time is $t = 25, 50$
  %  and $100$ respectively. The $x$-axis is number of moments $M$ and
  %  the $y$-axis is numerical damping rate. (c) The numerical
  %  frequencies with different number of moments with wave number
  %  $k=0.3$. The evolution time is $t = 25, 50$ and $100$
  %  respectively. The $x$-axis is number of moments $M$ and the
  %  $y$-axis is numerical frequency. The frequency of the electric
  %  field is estimated by counting the peaks in Figure
  %  \ref{fig:moment_03_a}.}
\end{figure}

\subparagraph{Moment convergence} 
Figure \ref{fig:moment_03_a} presents the time evolution of
$\mathcal{E}(t)$ of the filtered HME for different moment order
$M=40$, $50$, $60$, $70$ and $80$ with $k=0.3$, $N=3200$. One can
observe that the behavior of the time evolution of $\mathcal{E}(t)$
for different numbers of moments are all persisting on damping. This
indicates that the filter works well with different numbers of
moments.

The damping rates and frequencies for different end time $t_{end}$ are
also studied. Three end time $t_{end}=25$, $50$ and $100$ is studied,
and results are presented in Figure
\ref{fig:slope_frequency_03_moment_rate} and
\ref{fig:slope_frequency_03_moment_fre}.  It is clear that the damping
rates at different ending time are almost the same and they are
converging with the increasing of moment number.  Because of the
spatial error, there is a little distance between theoretical damping
rate and the converging limit. But as it is stated in
\ref{para:spatial}, the distance will go to zero with the increasing
of the spatial grid size.  Meanwhile, one can observe that the damping
frequencies for different ending time are just slightly
different. Moreover, with different numbers of moments, the
frequencies are almost same. This indicates that the filtered HME can
predict the correct damping rate and frequency with a small number of
moments.

\begin{figure}[!ht]
    \centering
    \psfrag{en}{\footnotesize{$\:\mathcal{E}$}}
    \psfrag{varivation}{\hspace{-12pt}relative variation}
    \includegraphics[width=.45\textwidth,height=5.0cm]{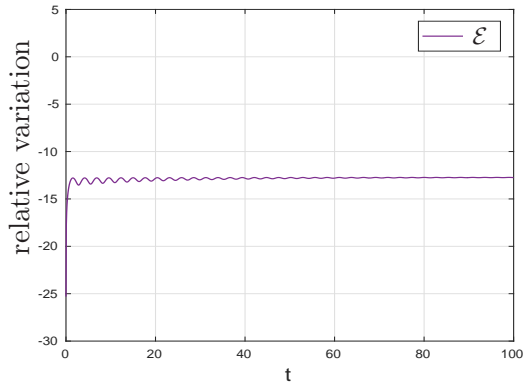}
    \caption{ \label{fig:conservative1d}
    Time evolution of the relative variation of the total
    energy $\mathcal{E}_{total}(t)$ in logarithm scale. The relative
    variation is defined by $|\mathcal{E}_{total}(t) -
    \mathcal{E}_{total}(0)| / \mathcal{E}_{total}(0)$.}
\end{figure}

\paragraph{Conservation property}
It has been proved in \cite{Wang} that the numerical scheme used in
this paper preserved the conservation of total mass and momentum,
expect the total energy.  We study the time evolution of the variation
of the total energy, which is presented in Figure
\ref{fig:conservative1d} for the setup $N=3200$, $M=50$ and $k=0.3$.
It is clear that the total energy changes very slightly in the whole
computation, which agrees with the result in \cite{Wang}. This
indicates the filter does not deteriorate the conservation of the
mass, momentum and energy.

\subsubsection{2D-2D example}
In this subsection, we study the filtered HME by 2D linear Landau
damping and try to claim that all the conclusions for 1D case also hold
for the 2D case. In the following, we pick up a few behaviors studied
for 1D case to validate our conclusion.
The initial value in \cite{Filbet_1} is adopted here:
\begin{equation}
  \label{eq:2d_initial}
  f(0,x,y,v_x,v_y) =f_{eq}=\frac{1}{2\pi}\exp\left(-\frac{v_x^2 +
      v_y^2}{2}\right)\Big(1 + A\cos(k_x x)\cos(k_y y)\Big),
\end{equation}
with $A= 10^{-3}$ and $k_x = k_y = 0.3$. The length of the periodic box
in the physical space is $L_x = L_y = 4 \pi / k_x$. The moment order
is set as $M= 40$. 

\begin{figure}[!ht]
  \centering
  \psfrag{E}{\footnotesize $\mathcal{E}$}
  \psfrag{en}{\footnotesize $\mathcal{E}$}
  \psfrag{variation}{\footnotesize\hspace{-14pt}relative variation}
  \subfigure[Filter effect]{
      \label{fig:2d_spatial_total_a}
    \begin{overpic}[ height=3.5cm]{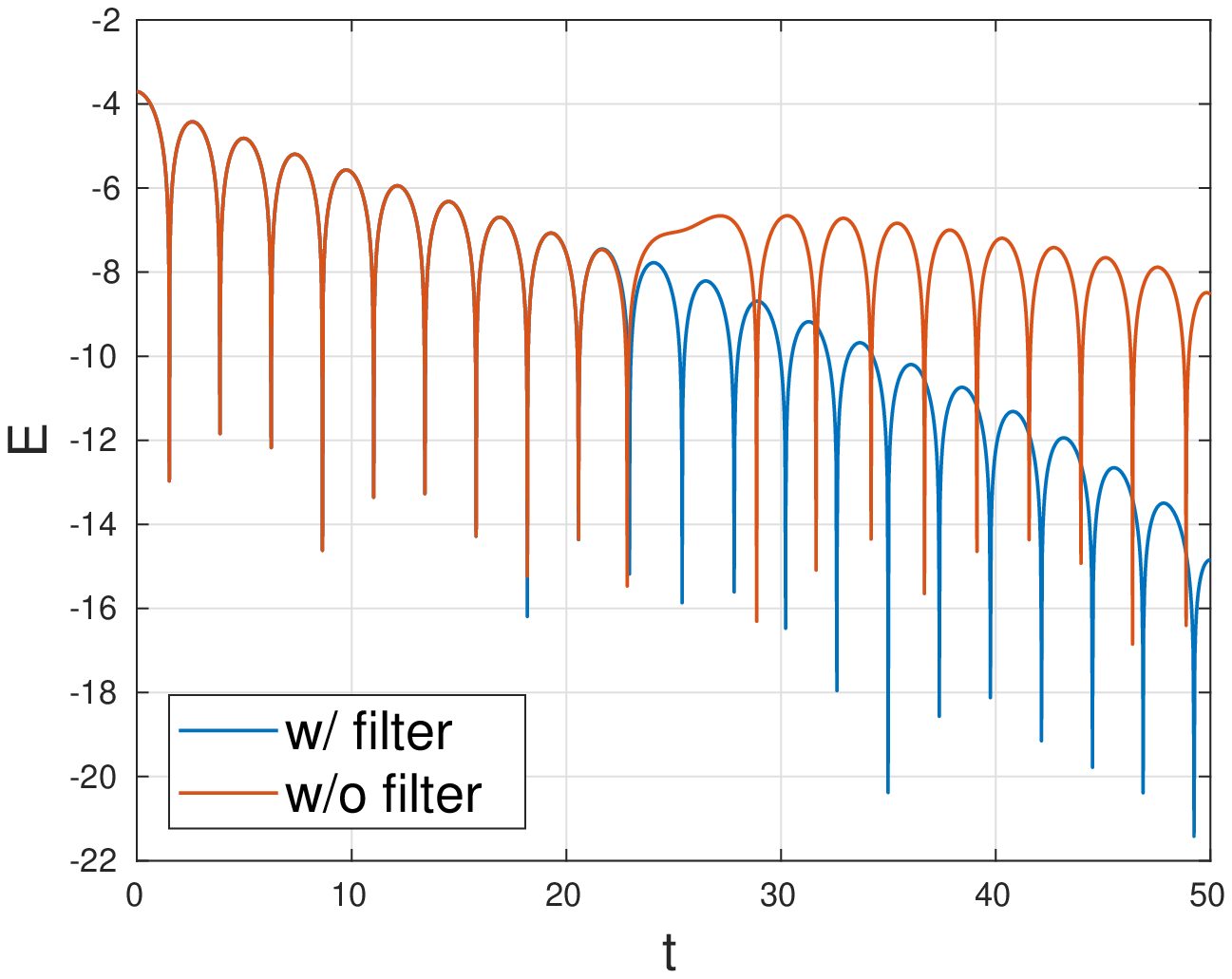}
    \end{overpic}
  } \quad 
  \subfigure[Landau damping]{
          \label{fig:2d_spatial_total_b}
    \begin{overpic}[height=3.5cm]{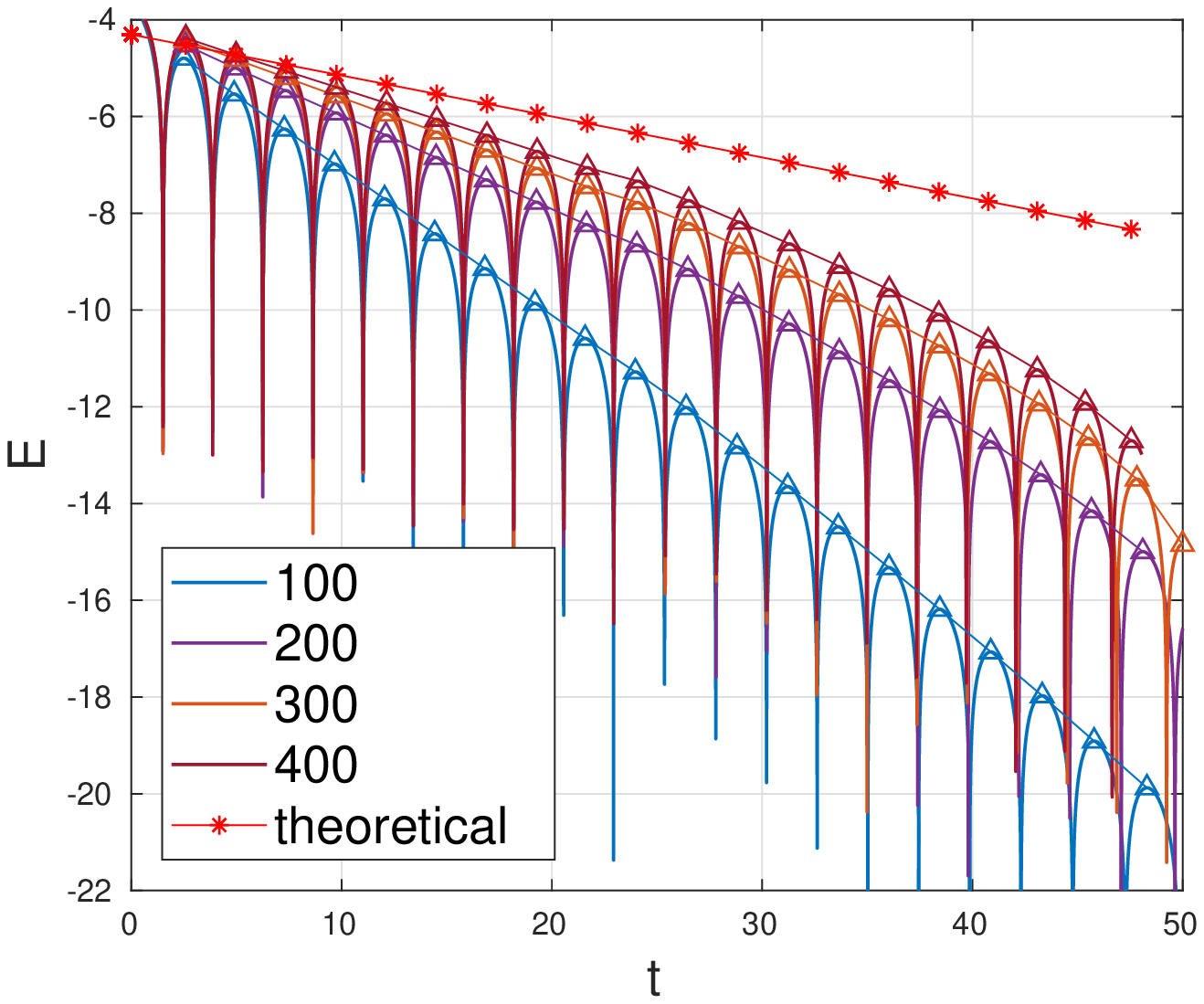}
    \end{overpic}
  }
  \subfigure[Conservation]{
    \label{fig:conservative_2d}
    \begin{overpic}[height=3.5cm]{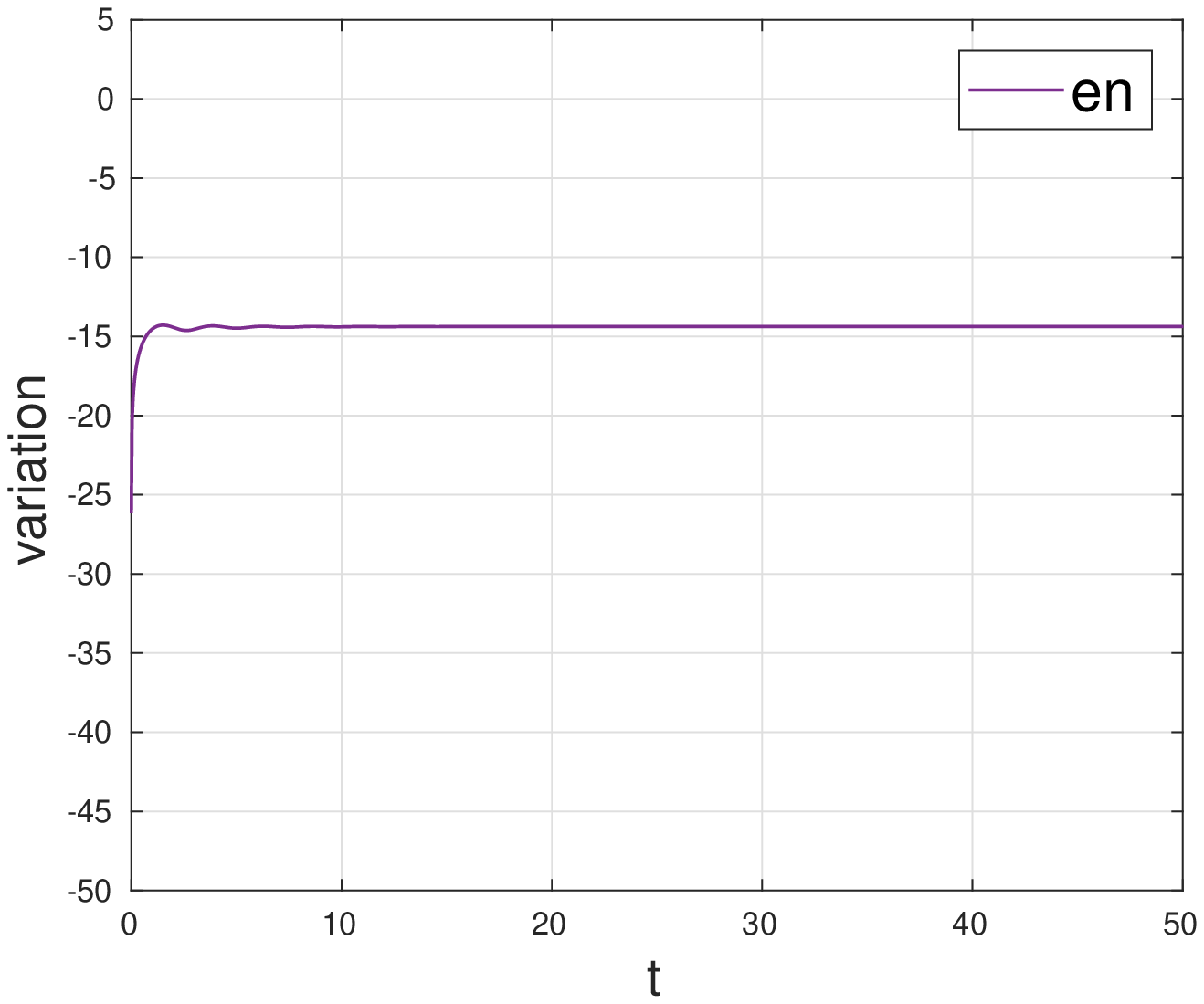}  
    \end{overpic}
  }
  \caption{\label{fig:2d_spatial_total} (a) Time evolution of
    $\ln(\mE(t))$ for HME and the filtered HME with $N_x=N_y=200$ and
    $M=40$.  (b) Time evolution of $\ln(\mE(t))$ for the filtered HME
    with $M=40$ for different grid sizes.  (c) Relative variation of
    the total energy.  }
  %(a) Exponential damping of $\mathcal{E}$ with the wave number
  %  $k = 0.3$. The blue line is that with filter while the red line is that
  %  without filter. Here $N_x=N_y=200$ and $M=40$. (b) Logarithm
  %  scale form of $\mathcal{E}$ in time with different spatial grid
  %  sizes. The slopes of the curves are the numerical damping rates by
  %  the least square fitting of the peak value points of
  %  $\mathcal{E}$. The slope of the star line is the damping rate given
  %  by the theoretic data.
\end{figure}

Figure. \ref{fig:2d_spatial_total_a} presents the time evolution of
the electric energy $\mE(t)$ of HME and the filtered HME. Just similar
to the 1D case, $\mE(t)$ for the filtered HME keeps on damping at the
time when the recurrence occurs for HME. Figure
\ref{fig:2d_spatial_total_b} shows the convergence of the damping rate
with the increasing of grid number, where the grid number are
$N = 100$, $200$, $300$ and $400$, respectively.

\begin{remark}
  The difference between the numerical result and the theoretical
  result may be caused by the smaller $N$ and not large enough moment
  order $M$, which is restricted by the computational cost.  Maybe the
  second order scheme we are working on can solve this problem.
\end{remark}

Figure \ref{fig:conservative_2d} illustrates the relative error of the
total energy with the time evolution for the 2D-2D case with the setup
$N_x= N_y = 100$, $M=40$ and $k=0.3$. Same as the 1D case, the total
energy is changing very slightly in the whole computation, which is
due to the splitting numerical method.

\subsection{Two-stream instability} 
In this subsection, the classical two-stream instability problem is
employed to study the filtered HME.
The two-stream instability is excited when the distribution function
is formed by two populations streaming in opposite directions with a
large enough relative drift velocity. Here, we adopt the same initial
condition as in
\cite{camporeale2016velocity}
\begin{equation}
  \label{eq:two_stream_ini}
  f(0, x, v) = \frac{1+ \epsilon \cos(kx)}{2\sqrt{2\pi}}
  \left(\exp\left(-\left(\frac{v+u_0}{\sqrt{2}u_{th}}\right)^2\right) 
    + \exp\left(-\left(\frac{v-u_0}{\sqrt{2}u_{th}}\right)^2\right)\right), 
\end{equation}
where the wave number $k = 0.5$, perturbation $\epsilon = 10^{-3}$,
velocity $u_0 = 1.0$, and thermal velocity $u_{th} = 0.5$
respectively.

\begin{figure}[!ht]
  \centering
  \psfrag{E}{\footnotesize $\mathcal{E}$}
  \subfigure[Spatial convergence]{
    \label{fig:two_stream_total}
    \begin{overpic}[height=5cm]{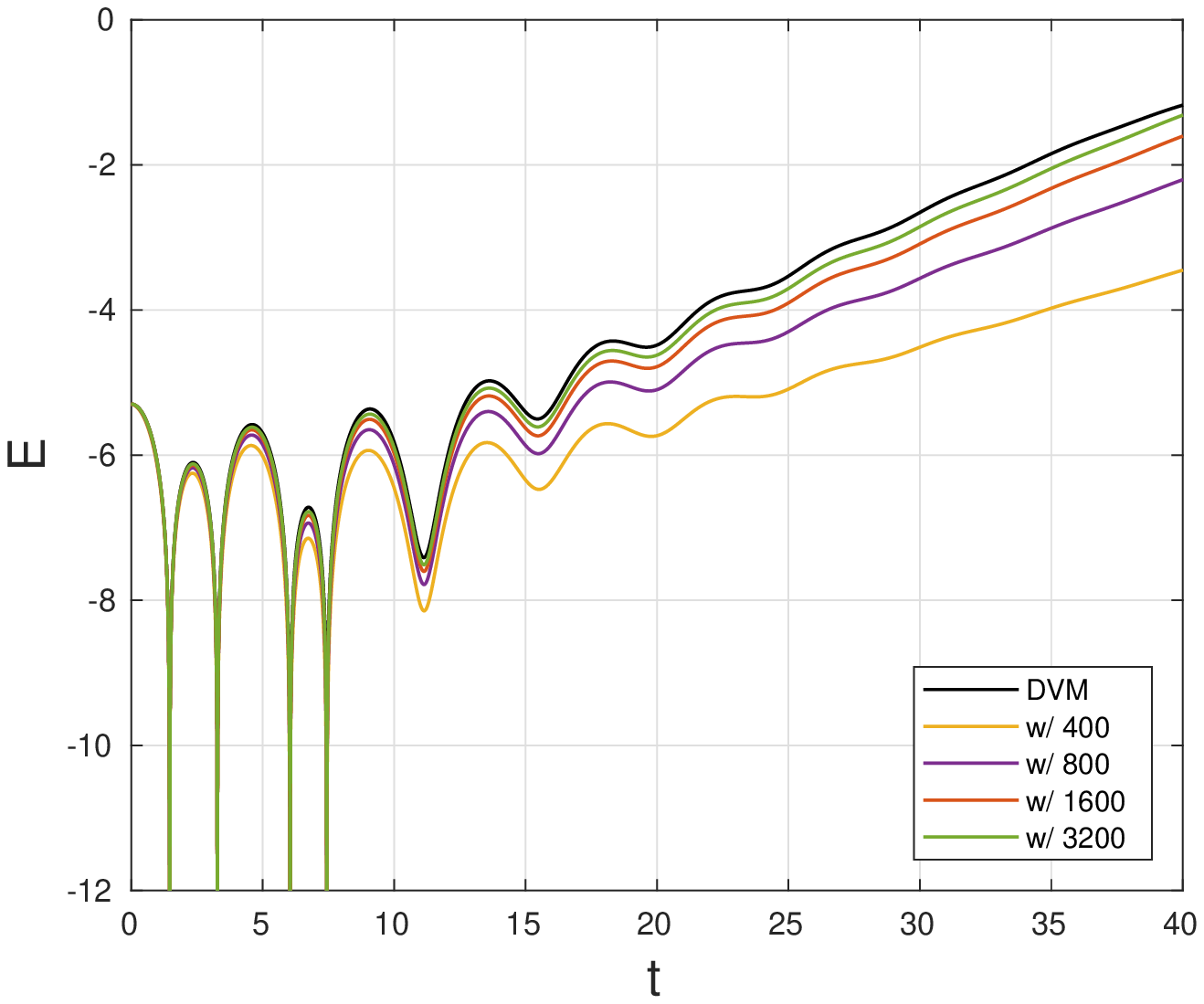}
    \end{overpic}
  }\quad
  \subfigure[With or without filter]{
    \label{fig:two_stream_local}
    \begin{overpic}[height=5cm]{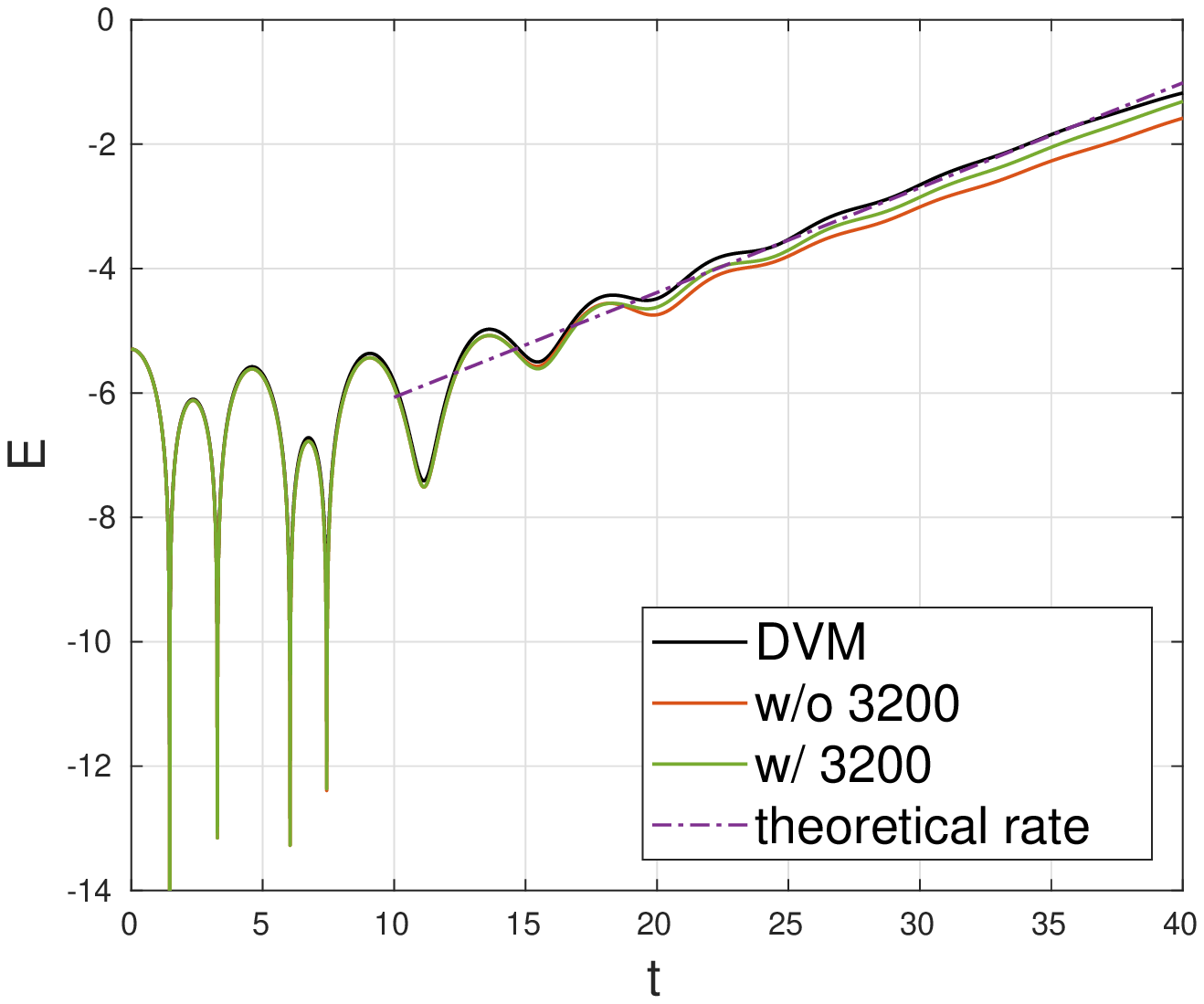}
    \end{overpic}
  }
  \caption{\label{fig:two_stream}
  (a) Time evolution of $\ln(\mE(t))$ of the filtered HME for
  different spatial grid sizes with $M=60$.
  (b) Time evolution of $\ln(\mE(t))$ of HME and the filtered HME 
  and comparison with reference solution.
  }
%    (a) The evolution of the electric energy
%    $\mathcal{E}$ in logarithm scale form by the moment method
%    without and with filter.  The black line is the reference
%    solution by the discrete velocity method. The dash lines are
%    those without filter, and the solid lines are those with
%    filter. Here the grid sizes are $N=400, 800, 1600, 3200$
%    respectively. The number of moment is $M=60$.  (b): The
%    comparisons of the numerical solutions without and with
%    filter. The black line is the reference solution by the discrete
%    velocity method. The dash line is that without filter, and the
%    solid line is that with filter. The dot-dash line is the
%    theoretical growth rate of the electric energy.  }
\end{figure}

Figure \ref{fig:two_stream} depicts the time evolution of the electric
energy $\mE(t)$, where the reference solution is computed by the
discrete velocity method (DVM) with the grid size large enough.
Figure \ref{fig:two_stream_total} shows the convergence of the
electric energy with respect to the grid size with $M=60$.  This
indicates the filtered HME can depict the two-stream instability
problem.  Figure \ref{fig:two_stream_local} presents the time
evolution of the electric energy of HME and the filtered HME and also
compare them with the reference solution. One can see the good
agreement with the solution of DVM and reference theoretical
rate. Moreover, the filtered HME behaviors better than HME for this
problem. Good agreement with the reference shows the power of the
filtered HME, saying this method is designed to suppress the numerical
recurrence, but it also works for the nonlinear two-steam instability
problem.

\begin{figure}[!ht]
  \centering
  \psfrag{en}{\footnotesize{$\:\mathcal{E}$}}
  \psfrag{variation}{\hspace{-12pt}relative variation}
  \includegraphics[width=.45\textwidth,height=5.0cm]{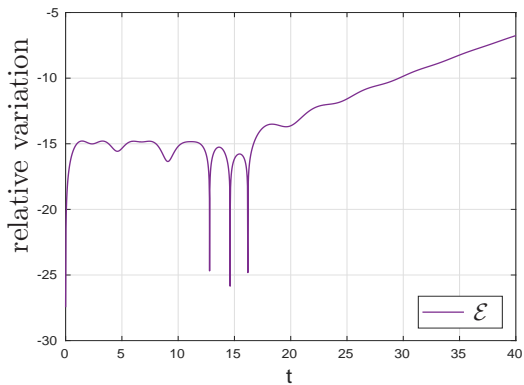}
  \caption{\label{fig:conservative_two-stream}
  Time evolution of the relative variation of the total energy $\mE_{total}(t)$ 
  for two-stream instability problem.
  }
  %The variation of the energy $\mathcal{E}_{total}(t)$ in
  %  time for the two-stream example.  The $x$-axis is time. The purple
  %  line corresponds to
  %  $\log(|\mathcal{E}_{total}(t) - \mathcal{E}_{total}(0)| /
  %  \mathcal{E}_{total}(0))$.}
\end{figure}

Figure \ref{fig:conservative_two-stream} illustrates the time
evolution of the relative error of the total energy. The setup is $N =
3200$ and $M=60$. It is shown that the total energy is changing very
slightly in the whole computation, and the variations of the total
energy in the nonlinear stage is also acceptable, since this is quite
a complicate physical process.

%%% Local Variables: 
%%% mode: latex
%%% TeX-master: "article"
%%% End: 

\section{Conclusion}\label{sec:conclusion}
We presented a filtered HME for the Vlasov-Poisson equations to
suppress the recurrence effects. Due to the careful construction, the
filter preserves most of physical properties of HME, including the
conservation of mass, momentum and energy, Galilean invariant and
convergence to VP. The quasi time-consistent property guarantees that
the solution to the filtered HME is not sensitive to the time step.
Two viewpoints on the quasi time-consistent filter show that the
filtered HME is a solver of the Vlasov equation, and it can predict
correct physical phenomena described by VP. Numerical simulations
demonstrate the power of the filter in suppressing recurrences and
producing more accurate solutions. The proposed numerical method can
depict the electric energy until it is close to the machine precision
for linear Landau damping cases and predict correct behaviors of the
electric energy for two-stream instability problem. More applications
of the filtered HME are in process.

\section*{Acknowledgements} 
This research of Y. Di is supported in part by the Natural Science
Foundation of China (Grant No. 11771437 and 91630208). And that of
Y. Wang is supported in part by the Natural Science Foundation of
China No. 11501042. R. Li is supported in part by the National Natural
Science Foundation of China (Grant No. 9163030002).

\bibliographystyle{plain}
\bibliography{article,tiao}

\begin{thebibliography}{10}

\bibitem{adjerid1986moving}
S.~Adjerid and J.~E. Flaherty.
\newblock A moving finite element method with error estimation and refinement
  for one-dimensional time dependent partial differential equations.
\newblock {\em SIAM J. Numer. Anal.}, 23(4):778--796, 1986.

\bibitem{armstrong1967}
T.~P. Armstrong.
\newblock Numerical studies of the nonlinear {V}lasov equation.
\newblock {\em Phys. Fluids}, 10:1269--1280, 1967.

\bibitem{armstrong1970solution}
T.~P. Armstrong, R.~C. Harding, G.~Knorr, and D.~Montgomery.
\newblock {S}olution of {V}lasov's equation by transform methods.
\newblock {\em J. Sci. Comput.}, 9:29--86, 1970.

\bibitem{BGK}
P.~L. Bhatnagar, E.~P. Gross, and M.~Krook.
\newblock A model for collision processes in gases. {I}. small amplitude
  processes in charged and neutral one-component systems.
\newblock {\em Phys. Rev.}, 94(3):511--525, 1954.

\bibitem{birdsall1985plasma}
Charles~K Birdsall and A~Bruce Langdon.
\newblock {\em Plasma physics via computer simulation}.
\newblock McGraw-Hill, NewYork, 2004.

\bibitem{Bourdiec2006Numerical}
S.~L. Bourdiec, F.~D. Vuyst, and L.~Jacquet.
\newblock Numerical solution of the vlasov–poisson system using generalized
  hermite functions.
\newblock {\em Commun. Comput. Phys.}, 175(8):528--544, 2006.

\bibitem{Fan_new}
Z.~Cai, Y.~Fan, and R.~Li.
\newblock Globally hyperbolic regularization of {G}rad's moment system.
\newblock {\em Comm. Pure Appl. Math.}, 67(3):464--518, 2014.

\bibitem{dvm2moments}
Z.~Cai, Y.~Fan, and R.~Li.
\newblock From discrete velocity model to moment method.
\newblock {\em Mathematica Numerica Sinica}, 38(3):227--244, 2016.

\bibitem{Tiao}
Z.~Cai, Y.~Fan, R.~Li, T.~Lu, and Y.~Wang.
\newblock Quantum hydrodynamic model by moment closure of wigner equation.
\newblock {\em J. Math. Phys.}, 53(10):103503, 2012.

\bibitem{Wang}
Z.~Cai, R.~Li, and Y.~Wang.
\newblock Solving {V}lasov equation using {NR$xx$} method.
\newblock {\em SIAM J. Sci. Comput.}, 35(6):A2807--A2831, 2013.

\bibitem{Cai2017AFilter}
Z.~Cai and Y.~Wang.
\newblock {S}uppression of recurrence in the {H}ermite-spectral method for
  transport equations.
\newblock {\em preprint}, 2017.

\bibitem{camporeale2016velocity}
E.~Camporeale, G.~L. Delzanno, B.~K. Bergen, and J.~D. Moulton.
\newblock On the velocity space discretization for the {V}lasov-{P}oisson
  system: {C}omparison between implicit {H}ermite spectral and
  {Particle-in-Cell} methods.
\newblock {\em Commun. Comput. Phys.}, 198:47--58, 2016.

\bibitem{canuto2012spectral}
C.~Canuto, M.~Y. Hussaini, A.~M. Quarteroni, A.~Thomas Jr, et~al.
\newblock {\em Spectral methods in fluid dynamics}.
\newblock Springer Science \& Business Media, 2012.

\bibitem{Carrillo2003}
J.~Carrillo, M.~Gamba, A.~Majorana, and C.~Shu.
\newblock A weno-solver for the transients of boltzmann-poisson system for
  semiconductor devices. performance and comparisons with monte carlo methods.
\newblock {\em J. Comput. Phys.}, 184:498--525, 2003.

\bibitem{Cheng}
C.~Z. Cheng and G.~Knorr.
\newblock The integration of the {V}lasov equation in configuration space.
\newblock {\em J. Comput. Phys.}, 22:330--351, 1976.

\bibitem{Cheng2013StudyOC}
Y.~Cheng, M.~Gamba, and J.~Morrison.
\newblock Study of conservation and recurrence of runge-kutta discontinuous
  galerkin schemes for vlasov-poisson systems.
\newblock {\em J. Sci. Comput.}, 56:319--349, 2013.

\bibitem{Filbet}
N.~Crouseilles and F.~Filbet.
\newblock Numerical approximation of collisional plasmas by high order methods.
\newblock {\em J. Comput. Phys.}, 201(2):546--572, 2004.

\bibitem{VM2015}
Y.~Di, Z.~Kou, and R.~Li.
\newblock High order moment closure for {V}lasov-{M}axwell equations.
\newblock {\em Front. Math. China}, 10(5):1087--1100, 2015.

\bibitem{eliasson2010numerical}
Bengt Eliasson.
\newblock Numerical simulations of the fourier-transformed vlasov-maxwell
  system in higher dimensions—theory and applications.
\newblock {\em Transport Theory and Statistical Physics}, 39(5-7):387--465,
  2010.

\bibitem{ellasson2001outflow}
B.~Ellasson.
\newblock {O}utflow boundary conditions for {F}ourier transformed
  one-dimensional {V}lasov-{P}oisson system.
\newblock {\em J. Sci. Comput.}, 16:1--28, 2001.

\bibitem{FATEMI1993209}
E.~Fatemi and F.~Odeh.
\newblock Upwind finite difference solution of boltzmann equation applied to
  electron transport in semiconductor devices.
\newblock {\em J. Comput. Phys.}, 108(2):209 -- 217, 1993.

\bibitem{Filbet2001Convergence}
F.~Filbet.
\newblock Convergence of a finite volume scheme for the vlasov--poisson system.
\newblock {\em SIAM J. Numer. Anal.}, 39(4):1146--1169, 2001.

\bibitem{Filbet2003Comparison}
F.~Filbet and E.~Sonnendr{\"u}cker.
\newblock Comparison of {E}ulerian {V}lasov solvers.
\newblock {\em Comput. Phys. Comm.}, 150(3):247--266, 2003.

\bibitem{Filbet_1}
F.~Filbet, E.~Sonnendr{\"u}cker, and P.~Bertrand.
\newblock Conservative numerical schemes for the {V}lasov equation.
\newblock {\em J. Comput. Phys.}, 172:166--187, 2001.

\bibitem{Gottlieb2001spectral}
D.~Gottlieb and J.~S. Hesthaven.
\newblock Spectral methods for hyperbolic problems.
\newblock {\em J. Comput. Appl. Math.}, 128:83--131, 2001.

\bibitem{Grad}
H.~Grad.
\newblock On the kinetic theory of rarefied gases.
\newblock {\em Comm. Pure Appl. Math.}, 2(4):331--407, 1949.

\bibitem{grant1967fourier}
F.~C. Grant and M.~R. Feix.
\newblock {F}ourier-{H}ermite solutions of the {V}lasov equations in the
  linearized limit.
\newblock {\em Phy. Fluids}, 10(4):696--702, 1967.

\bibitem{heath2012discontinuous}
R.~E. Heath, I.~M. Gamba, P.~J. Morrison, and C.~Michler.
\newblock A discontinuous {G}alerkin method for the {V}lasov-{P}oisson system.
\newblock {\em J. Comput. Phys.}, 231(4):1140--1174, 2012.

\bibitem{Hesthaven2008filtering}
J.~S. Hesthaven and R.~Kirby.
\newblock {F}iltering in {L}egendre spectral methods.
\newblock {\em Math. Comput.}, 77(263):1425--1452, 2008.

\bibitem{Holloway1996Spectral}
J.~P. Holloway.
\newblock Spectral velocity discretizations for the {V}lasov-{M}axwell
  equations.
\newblock {\em Transport. Theor. Stat.}, 25(1):1--32, 1996.

\bibitem{HouLi2007}
T.~Hou and R.~Li.
\newblock Computing nearly singular solutions using pseudo-spectral methods.
\newblock {\em J. Comput. Phys.}, 226(1):379--397, 2007.

\bibitem{Joyce1971Numerical}
G.~Joyce, G.~Knorr, and H.~K. Meier.
\newblock Numerical integration methods of the {Vlasov} equation.
\newblock {\em J. Comput. Phys.}, 8(1):53--63, 1971.

\bibitem{Kanevsky2006Idempotent}
A.~Kanevsky, K.~Carpenter, and J.~S. Hesthaven.
\newblock {I}dempotent filtering in spectral and spectral element methods.
\newblock {\em J. Comput. Phys}, 220(1):41 -- 58, 2006.

\bibitem{Klimas1987method}
A.~J. Klimas.
\newblock A method for overcoming the velocity space filamentation problem in
  collisionless plasma model solutions.
\newblock {\em J. Comput. Phys.}, 68(1):202--226, 1987.

\bibitem{Klimas1994splitting}
A.~J. Klimas and W.~M. Farrell.
\newblock A splitting algorithm for {V}lasov simulation with filamentation
  filtration.
\newblock {\em J. Comput. Phys.}, 110(1):150--163, 1994.

\bibitem{Kreiss1979Stability}
H.O. Kreiss and J.~Oliger.
\newblock {S}tability of the {F}ourier method.
\newblock {\em SIAM J. Numer. Anal.}, 16:421 -- 433, 1979.

\bibitem{Landau1946On}
L.~Landau.
\newblock {O}n the vibrations of the electronic plasma.
\newblock {\em Eur. J. Org. Chem.}, 2006(2):498–506, 1946.

\bibitem{mcclarren2010robust}
R.~G. McClarren and C.~D. Hauck.
\newblock Robust and accurate filtered spherical harmonics expansions for
  radiative transfer.
\newblock {\em J. Comput. Phys.}, 229(16):5597--5614, 2010.

\bibitem{Muller}
I.~M{\"u}ller and T.~Ruggeri.
\newblock {\em Rational Extended Thermodynamics, Second Edition}, volume~37 of
  {\em Springer tracts in natural philosophy}.
\newblock Springer-Verlag, New York, 1998.

\bibitem{ng2004complete}
C.~S. Ng, A.~Bhattacharjee, and F.~Skiff.
\newblock Complete spectrum of kinetic eigenmodes for plasma oscillations in a
  weakly collisional plasma.
\newblock {\em Phys. Rev. Lett.}, 92(6):065002, 2004.

\bibitem{parker2015fourier}
J.~T. Parker and P.~J. Dellar.
\newblock {F}ourier-{H}ermite spectral representation for the
  {V}lasov-{P}oisson system in the weakly collisional limit.
\newblock {\em J. Plasma. Phys.}, 81(02):305810203, 2015.

\bibitem{Qiu2011positivity}
J.~Qiu and C.~Shu.
\newblock {P}ositivity preserving semi-{L}agrangian discontinuous {G}alerkin
  formulation: {T}heoretical analysis and application to the {V}lasov-{P}oisson
  system.
\newblock {\em J. Comput. Phys.}, 230(23):8386 -- 8409, 2011.

\bibitem{UTH}
J.~W. Schumer and J.~P. Holloway.
\newblock Vlasov simulation using velocity-scaled {H}ermite representations.
\newblock {\em J. Comput. Phys.}, 144(2):626--661, 1998.

\bibitem{Shoucri1974Numerical}
M.~Shoucri and G.~Knorr.
\newblock Numerical integration of the vlasov equation ☆.
\newblock {\em J. Comput. Phys.}, 14(1):84--92, 1974.

\bibitem{SL}
E.~Sonnendr{\"u}cker, J.~Roche, P.~Betrand, and A.~Ghizzo.
\newblock The semi-{L}agrangian method for the numerical resolution of {V}lasov
  equations.
\newblock {\em J.Comput. Phys}, 149(2):201--220, 1998.

\bibitem{Torrilhon2006}
M.~Torrilhon.
\newblock Two dimensional bulk microflow simulations based on regularized
  {G}rad's 13-moment equations.
\newblock {\em SIAM Multiscale. Model. Simul.}, 5(3):695--728, 2006.

\bibitem{Vlasov}
A.~A. Vlasov.
\newblock On vibration properties of electron gas.
\newblock {\em J. Exp. Theor. Phys.}, 8(3):291, 1938.

\bibitem{zaki1988finiteI}
S.~I. Zaki, R.~T. Gardner, and T.~J. Boyd.
\newblock {A} finite element code for the simulation of one-dimensional
  {V}lasov plasmas. i. {T}heory.
\newblock {\em J. Comput. Phys.}, 79:184--199, 1988.

\end{thebibliography}
\end{document}